\documentclass{amsart}

\usepackage{pictexwd}

\theoremstyle{plain}
\newtheorem{theorem}{Theorem}
\newtheorem{lemma}{Lemma}

\newtheorem{remark}{Remark}
\newtheorem{definition}{Definition}

\begin{document}

\title{On Hopf's Lemma and the Strong Maximum Principle}

\author{S. Bertone}
\address{Dipartimento di Matematica e Applicazioni,
         Universit\'{a} degli Studi di Milano-Bicocca,
         Via R. Cozzi 53,
         20125 Milano}
\email{bertone@matapp.unimib.it}

\author{A. Cellina}
\address{Dipartimento di Matematica e Applicazioni,
         Universit\'{a} degli Studi di Milano-Bicocca,
         Via R. Cozzi 53,
         20125 Milano}
\email{cellina@matapp.unimib.it}

\author{E. M. Marchini}
\address{Dipartimento di Matematica e Applicazioni,
         Universit\'{a} degli Studi di Milano-Bicocca,
         Via R. Cozzi 53,
         20125 Milano}
\email{elsa@matapp.unimib.it}


\date{}


\begin{abstract}
In this paper we consider Hopf's Lemma and the Strong Maximum Principle for supersolutions to
$$\sum_{i=1}^N g_i(u_{x_i}^2)u_{x_ix_i}=0$$
under suitable hypotheses that allow $g_i$ to assume value zero at zero.
\end{abstract}

\maketitle

\section{Introduction}

Let $\Omega\subset\mathbb R^N$ be a connected, open and bounded set; we call $\Omega$ regular if for every $z\in\partial\Omega$, there exists a tangent plane, continously depending on $z$. We say that $\Omega $ satisfies the interior ball condition at $z$ if there exists an open ball $B\subset\Omega$ with $z\in\partial B$.
On $\Omega$, consider the operator
\begin{equation}\label{elJ}
F(u)=\sum_{i=1}^N g_i(u_{x_i}^2)u_{x_ix_i},
\end{equation}
where $g_i:[0,+\infty)\to[0,+\infty)$ are continuous functions.

When $F$ is elliptic, two classical results hold.\\

{\it Hopf's Lemma:\\
Let $\Omega$ be regular, let $u$ be such that $F(u)\leq0$ on $\Omega$. Suppose that there exists $z\in\partial\Omega$ such that
$$u(z)<u(x),\quad\mbox{ for all $x$ in $\Omega$}.$$
If, in addition, $\Omega$ satisfies the interior ball condition at $z$, we have
$$\frac{\partial u}{\partial\nu}(z)<0,$$
where $\nu$ is the outer unit normal to $B$ at $z$.\\}

{\it The Strong Maximum Principle:\\
Let $u$ be such that $F(u)\leq0$ on $\Omega$, then if $u$ it attains minimum in $\Omega$, it is a constant.}\\

In 1927 Hopf proved the Strong Maximum Principle in the case of second order elliptic partial differential equations, by applying a comparison technique, see \cite{h}. For the class of quasilinear elliptic problems, many contributions have been given, to extend the validity of the previous results, as in
\cite{bnv,dpr,d,fq,gt,g,ps,ps2,psz,se,se2,v}.

In the case in equation (\ref{elJ}) we have $g_i\equiv1$, for every $i$, then $F(u)=\Delta u$, and we find the classical problem of the Laplacian, see \cite{e,gt}.

On the other hand, when there exists $i\in\{1,\dots,N\}$ such that $g_i\equiv0$ on an interval $I=[0,T]\subset\mathbb R$, the Strong Maximum Principle does not hold. Indeed, in this case, it is always possible to define a function $u$ assuming minimum in $\Omega$ and such that $\sum_{i=1}^N g_i(u_{x_i}^2)u_{x_ix_i}=0$. For instance, let $g_N(t)=0$, for every $t\in[0,2]$. The function
$$u(x_1,\dots,x_N)=\left\{\begin{array}{ll}
        -(x_N^2-1)^4 & \mbox{ if }-1\leq x_N\leq1\\
        0 & \mbox{ otherwise}
        \end{array}
\right.,$$
satisfies (\ref{elJ}) in $\mathbb R^N$.

We are interested in the case when $0\leq g_i(t)\leq1$ and it does not exist $i$ such that $g_i\equiv0$ on an interval.
Since $g_i$ could assume value zero, the equation (\ref{elJ}) is non elliptic.

The results known so far, for the validity of Hopf's Lemma and of the Strong Maximum Principle, suggest that, for possibly non elliptic equations, but arising from a functional having rotational symmetry,  this validity shall depend only on the behaviour of the functions $g_i$ near zero, see \cite{cel}.

In this paper, we prove, in section 3, a sufficient condition for the validity of the Hopf's Lemma and of the Strong Maximum Principle; a remarkable feature of this condition is that it concerns only the behaviour of the function $g_i(t)$ that goes fastest to zero, as $t$ goes to zero.
Hopf's lemma and the Strong Maximum Principle are essentially the same result as long as we can build subsolutions whose level lines can have arbitrarily large curvature. This need not be always possible for problems not possessing rotational symmetry. This difficulty will be evident in sections 4 and 5.
In these sections, a more restricted class of equations is considered, namely when all the functions $g_i$, for $i=1,\dots,N-1$, are $1$ and only $g_N$ is allowed to go to zero.
In this simpler class of equations we are able to show that the condition
$$\lim_{t\rightarrow0^+}\frac{\left(g_N(t)\right)^{3/2}}{tg_N^\prime(t)}>0$$
is at once necessary for the validity of Hopf's Lemma and sufficient for the validity of the Strong Maximum Principle.

\section{Preliminary results}

We impose the following local assumptions.\\

{\it Assumptions (L):\\
\\
There exists $\overline t>0$ such that:\\

i) on $[0,\overline t\:]$, for every $i=1,\dots,N-1$,
$$0\leq g_N(t)\leq g_i(t)\leq1;$$

ii) $g_N$ is continuous on $[0,\overline t\,]$; positive and differentiable on $(0,\overline t\,]$;\\

iii) on $(0,\overline t\,]$, the function $t\to g_N(t)+g_N^\prime(t)t$ is non decreasing.}\\

Notice that, in case {\it ii)} above is violated, the Strong Maximum Principle does not hold; and that condition {\it iii)} above includes the case of the Laplacian, $g_i(t)\equiv1$; and, finally, that under these assumptions, $g_i$ could assume value zero at most for $t=0$.

Moreover, we can consider the equation
$$\sum_{i=1}^N g_i(u_{x_i}^2)u_{x_ix_i}=0$$
as the Euler-Lagrange equation associated to the functional
$$J(u)=\int_\Omega{L(\nabla u)d\Omega}=\int_\Omega{\frac{1}{2}\left(\sum_{i=1}^N f_i(u_{x_i}^2)u_{x_i}^2\right)d\Omega},$$
where $L(\nabla u)$ is strictly convex in $\{(u_{x_1},\dots,u_{x_N}):u_{x_i}^2\leq\overline t, \mbox{ for every }i=1,\dots,N\}$. Indeed, fix $i$. Let $f_i$ be a solution to the differential equation
\begin{equation}\label{eqg}
g_i(t)=f_i(t)+5tf_i^\prime(t)+2t^2f_i^{\prime\prime}(t),
\end{equation}
for $t\in[0,\overline t\,]$.
Since
$$\frac{\partial^2 L}{\partial u_{x_i}^2}(u_{x_i}^2)=f_i(u_{x_i}^2)+5u_{x_i}^2f_i^\prime(u_{x_i}^2)+2u_{x_i}^4f_i^{\prime\prime}(u_{x_i}^2)=g_i(u_{x_i}^2),$$
we have that
$$\mbox{div}\nabla_{\nabla u}L(\nabla u)=\sum_{i=1}^N\frac{\partial^2 L}{\partial u_{x_i}^2}u_{x_ix_i}=\sum_{i=1}^N g_i(u_{x_i}^2)u_{x_ix_i}.$$
Moreover, the strict convexity of $L(\nabla u)$ in $\{(u_{x_1},\dots,u_{x_N}):u_{x_i}^2\leq\overline t,\mbox{ for every }i=1,\dots,N\}$ follows by the fact that $g_i$ is positive in $(0,\overline t\,]$.\\

Since we will need general comparison theorems that depend on the global properties of the solutions, i.e. on their belonging to a Sobolev space, we will need also a growth assumption on $g_i$ (assumption {\it (G)}) to insure these properties of the solutions.\\

{\it Assumption (G):\\
Each function $f_i$ as defined in (\ref{eqg}), is bounded and $f_i(u_{x_i}^2)u_{x_i}^2$ is strictly convex.
}\\

Any function $g_i$ satisfying assumptions {\it (L)} on $[0,\overline t\,]$
can be extended so as to satisfy assumption {\it (G)} on $[0,+\infty)$. In fact, it is enough to extend $g_i$ to $(\,\overline t,+\infty)$ by setting $g_i(t)=f_i(\,\overline t\,)$, for $t>\overline t$.\\

\begin{definition}

Let $\Omega$ be open, and let $u\in W^{1,2}(\Omega)$. The map $u$ is a weak solution to the equation $F(u)=0$ if, for every $\eta\in C^{\infty}_0(\Omega)$,
$$\int_\Omega{\langle\nabla L(\nabla u(x)),\nabla\eta(x)\rangle dx}=0.$$
$u$ is a weak subsolution ($F(u)\geq0$) if, for every $\eta\in C^{\infty}_0(\Omega)$, $\eta\geq0$,
$$\int_\Omega{\langle\nabla L(\nabla u(x)),\nabla\eta(x)\rangle dx}\leq0.$$
$u$ is a weak supersolution ($F(u)\leq0$) if, for every $\eta\in C^{\infty}_0(\Omega)$, $\eta\geq0$,
$$\int_\Omega{\langle\nabla L(\nabla u(x)),\nabla\eta(x)\rangle dx}\geq0.$$
We say that a function $w\in W^{1,2}(\Omega)$ is such that $w_{|\partial \Omega}\leq0$ if $w^{+}\in W^{1,2}_0(\Omega)$.

\end{definition}

The growth assumption {\it (G)} assures that, if $u\in W^{1,2}(\Omega)$, then $\nabla L(\nabla u(x))\in L^2(\Omega)$. The strict convexity of $L$ implies the following comparison lemma.

\begin{lemma}\label{comp}

Let $\Omega$ be a open and bounded set, let $v\in W^{1,2}(\Omega)$ be a subsolution and let $u\in W^{1,2}(\Omega)$ be a supersolution to the equation $F(u)=0$.
If $v_{|\partial \Omega}\leq u_{|\partial \Omega}$, then $v\leq u$ a.e. in $\Omega$.

\end{lemma}

We wish to express the operator
$$F(v)=\sum_{i=1}^N g_i(v_{x_i}^2)v_{x_ix_i}$$
in polar coordinates. Set
$$\left\{\begin{array}{l}
        x_1=\rho\cos\theta_{N-1}\dots\cos\theta_2\cos\theta_1\\
        x_2=\rho\cos\theta_{N-1}\dots\cos\theta_2\sin\theta_1\\
        \dots\\
        x_N=\rho\sin\theta_{N-1}
        \end{array}
\right.$$
so that
$$v_{x_i}=v_{\rho}\frac{x_i}{\rho}\quad\mbox{ and }\quad v_{x_ix_i}=v_{\rho\rho}\left(\frac{x_i}{\rho}\right)^2+\frac{v_{\rho}}{\rho}\left[1-\left(\frac{x_i}{\rho}\right)^2\right].$$
When $v$ is a radial function, $F$ reduces to
$$F(v)=\sum_{i=1}^Ng_i\left(v_{\rho}^2\left(\frac{x_i}{\rho}\right)^2\right)\left[v_{\rho\rho}\left(\frac{x_{i}}{\rho}\right)^{2}+\frac{v_{\rho}}{\rho}\left(1-\left(\frac{x_{i}}{\rho}\right)^{2}\right)\right]=$$
$$v_{\rho\rho}\sum_{i=1}^Ng_i\left(v_{\rho}^2\left(\frac{x_i}{\rho}\right)^2\right)\left(\frac{x_i}{\rho}\right)^2+\frac{v_{\rho}}{\rho}\sum_{i=1}^Ng_i\left(v_{\rho}^2\left(\frac{x_i}{\rho}\right)^2\right)\left(1-\left(\frac{x_i}{\rho}\right)^2\right).$$
In general, we do not expect that the equation $F(v)=0$ admits radial solutions. However we will use the expression of $F$ valid for radial functions in order to reach our results.\\

The following technical lemmas will be used later.

\begin{lemma}\label{ln}

Let $n=2,\dots,N$ and set
$$h_n(a)=g_N\left(\frac{t(1-a)}{n-1}\right)(1-a)+g_N(ta)a.$$
For every $0<t\leq\overline t$ (\;$\overline t$ defined in assumptions (L)), $h_n(a)\geq h_n(1/n)$, for every $a$ in $[0,1]$.

\end{lemma}

\proof

Since, on $(0,\overline t\,]$, the function $t\to g_N(t)+g_N^\prime(t)t$ is non decreasing, we have that
$$h_n^{\prime}(a)=-g_N\left(\frac{t(1-a)}{n-1}\right)-g^\prime_N\left(\frac{t(1-a)}{n-1}\right)\frac{t(1-a)}{n-1}+g_N(ta)+g_N^\prime(ta)ta\geq0$$
if and only if $a\geq1/n$, so that $h_n(a)\geq h_n(1/n)$, for every $a\in[0,1]$.

\endproof

\begin{lemma}\label{lind}

For every $0<t\leq\overline t$ (\;$\overline t$ defined in assumptions (L)), we have that
$$\sum_{i=1}^N g_N\left(t\left(\frac{x_i}{\rho}\right)^2\right)\left(\frac{x_i}{\rho}\right)^2\geq g_N\left(\frac{t}{N}\right).$$

\end{lemma}

\proof

We prove the claim by induction on $N$.

Let $N=2$. Set $a=\sin^2\theta_1$. Applying Lemma \ref{ln} we obtain that
$$g_N\left(t\left(\frac{x_1}{\rho}\right)^2\right)\left(\frac{x_1}{\rho}\right)^2+g_N\left(t\left(\frac{x_2}{\rho}\right)^2\right)\left(\frac{x_2}{\rho}\right)^2=$$
$$g_N(t(1-a))(1-a)+g_N(ta)a\geq g_N\left(\frac{t}{2}\right).$$

Suppose that the claim is true for $N-1$, i.e.
$$\sum_{i=1}^{N-1}g_N\left(t\left(\frac{x_i}{\rho}\right)^2\right)\left(\frac{x_i}{\rho}\right)^2\geq g_N\left(\frac{t}{N-1}\right).$$
Let us prove it for $N$. Set
$$\left\{\begin{array}{l}
        y_1=\rho\cos\theta_{N-2}\dots\cos\theta_2\cos\theta_1\\
        y_2=\rho\cos\theta_{N-2}\dots\cos\theta_2\sin\theta_1\\
        \dots\\
        y_{N-1}=\rho\sin\theta_{N-2}
        \end{array}
\right.$$
and set
$a=\sin^2\theta_{N-1}$. Applying Lemma \ref{ln} we obtain that
$$\sum_{i=1}^{N}g_N\left(t\left(\frac{x_i}{\rho}\right)^2\right)\left(\frac{x_i}{\rho}\right)^2=$$
$$\sum_{i=1}^{N-1}g_N\left(t\left(\frac{y_i}{\rho}\right)^2(1-a)\right)\left(\frac{y_i}{\rho}\right)^2(1-a)+g_N(ta)a\geq$$ $$g_N\left(\frac{t(1-a)}{N-1}\right)(1-a)+g_N(ta)a\geq g_N\left(\frac{t}{N}\right),$$
and the claim is proved.

\endproof

\section{A sufficient condition for the validity of Hopf's Lemma and of the Strong Maximum Principle}

Consider the improper Riemann integral
$$\int^\xi_0{\frac{g_N(\zeta^2/N)}{\zeta}d\zeta}=\lim_{\widehat{\xi}\to 0}\int^\xi_{\widehat{\xi}}{\frac{g_N(\zeta^2/N)}{\zeta}d\zeta}$$
as an extended valued function $G$,
$$G(\xi)=\int^\xi_0{\frac{g_N(\zeta^2/N)}{\zeta}d\zeta},$$
where we mean that $G(\xi)\equiv+\infty$ whenever the integral diverges.\\

We wish to prove the following lemma.

\begin{lemma}[Hopf's Lemma]\label{hopfr}

Let $\Omega\subset\mathbb R^N$ be a connected, open and bounded set. Let $u\in W^{1,2}(\Omega)\cap C\left(\overline\Omega\right)$ be a weak solution to
$$\sum_{i=1}^N g_i(u_{x_i}^2)u_{x_ix_i}\leq0.$$
In addition to the assumptions (L) and (G) on $g_i$, assume that $G(\xi)\equiv+\infty$.
Suppose that there exists $z\in\partial\Omega$ such that
$$u(z)<u(x),\quad\mbox{ for all $x$ in $\Omega$}$$
and that $\Omega$ satisfies the interior ball condition at $z$. Then
$$\frac{\partial u}{\partial\nu}(z)<0,$$
where $\nu$ is the outer unit normal to $B$ at $z$.

\end{lemma}

As an example of an equation satisfying the assumptions of the theorem above, consider the Laplace equation $\Delta u=0$. The functions $g_i\equiv1$ satisfy the assumptions {\it (L)} and {\it (G)}, and
$$G(\xi)=\int^\xi_0{\frac{1}{\zeta}d\zeta}=+\infty.$$

Another example is obtained setting
$$g_N(t)=\dfrac{1}{|\ln(t)|}$$
for $0\leq t\leq1/e$. The assumptions {\it (L)} and {\it (G)} are satisfied; moreover, for $0\leq\xi^2/N\leq1/e$,
$$G(\xi)=\int^{\xi}_{0}{\frac{d\zeta}{\zeta|\ln(\zeta^2/N)|}}=+\infty.$$\\

\proof[Proof of Lemma \ref{hopfr}]

a) Assume that $u(z)=0$ and that $B=B(O,r)$. We prove the claim by contradiction. Suppose that
$$\frac{\partial u}{\partial\nu}(z)\geq0,$$
where $\nu$ is the outer unit normal to $B$ at $z$.
Let $\epsilon=\min\left\{u(x):x\in\overline{B(O,r/2)}\right\}$; we have that $\epsilon>0$. Set
$$\omega=B(O,r)\setminus\overline{B(O,r/2)}.$$

b) We seek a radial function $v\in W^{1,2}(\omega)\cap C(\overline{\omega})$ satisfying
\begin{equation}\label{eqnrad}
\left\{\begin{array}{ll}
        v \mbox{ is a weak solution to }F(v)\geq0 & \mbox{ in } \omega\\
        v>0 & \mbox{ in }\omega\\
        v=0 & \mbox{ in }\partial B(O,r)\\
        v\leq\epsilon & \mbox{ in }\partial B(O,r/2)\\
        v_\rho(z)<0.\\
        \end{array}
\right.
\end{equation}
Consider the Cauchy problem
\begin{equation}\label{cr}
\left\{\begin{array}{l}
        \zeta^\prime=-\dfrac{N-1}{\rho}\dfrac{\zeta}{g_N(\zeta^2/N)}\\
        \\
        \zeta(r/2)=-\dfrac{\epsilon}{r}.\\
\end{array}
\right.
\end{equation}
There exists a unique local solution $\zeta$ of (\ref{cr}), such that
$$\int^{\zeta(\rho)}_{-\frac{\epsilon}{r}}{\frac{g_N(\zeta^2/N)}{\zeta}d\zeta}=\int^{\rho}_{r/2}{-\frac{N-1}{s}ds}=-(N-1)\ln\left(\frac{2\rho}{r}\right).$$
We claim that $\zeta$ is defined in $[r/2,+\infty)$. Indeed, suppose that $\zeta$ is defined in $[r/2,\tau)$, with $\tau<+\infty$. Since $\zeta^\prime>0$, $\zeta$ is an increasing function, so that $\tau<+\infty$ if and only if $\lim_{\rho\to\tau}\zeta(\rho)=0$. But
$$-\infty=\lim_{\rho\to\tau}\int^{\zeta(\rho)}_{-\frac{\epsilon}{r}}{\frac{g_N(\zeta^2/N)}{\zeta}d\zeta}=\lim_{\rho\to\tau}-(N-1)\ln\left(\frac{2\rho}{r}\right),$$
a contradiction. Hence, the solution $\zeta$ of (\ref{cr}) is defined in $[r/2,+\infty)$.

Setting $v_\rho=\zeta$, since, for every $\rho\in(r/2,r)$,
$$-\frac{\epsilon}{r}<v_\rho(\rho)<0,$$
we have that the function
$$v(\rho)=\int_{r}^\rho{v_{\rho}(s)ds}$$
solves the problem
$$v_{\rho\rho}g_N\left(\frac{v_{\rho}^2}{N}\right)+\frac{v_{\rho}}{\rho}(N-1)=0,$$
in particular, $v(\rho)>0$ and $v_\rho(\rho)<0$, for every $\rho\in(r/2,r)$, $v(r)=0$ and $v(r/2)\leq\epsilon$.

Since $v_{\rho\rho}\geq0$ and $-\sqrt{\overline{t}}\leq v_{\rho}\leq0$, for every $\rho\in(r/2,r)$, by the hypotheses on $g_i$ and by Lemma \ref{lind}, we have that
$$F(v)=v_{\rho\rho}\sum_{i=1}^Ng_i\left(v_{\rho}^2\left(\frac{x_i}{\rho}\right)^2\right)\left(\frac{x_i}{\rho}\right)^2+\frac{v_{\rho}}{\rho}\sum_{i=1}^Ng_i\left(v_{\rho}^2\left(\frac{x_i}{\rho}\right)^2\right)\left(1-\left(\frac{x_i}{\rho}\right)^2\right)\geq$$ $$v_{\rho\rho}\sum_{i=1}^Ng_N\left(v_{\rho}^2\left(\frac{x_i}{\rho}\right)^2\right)\left(\frac{x_i}{\rho}\right)^2+\frac{v_{\rho}}{\rho}\sum_{i=1}^N\left(1-\left(\frac{x_i}{\rho}\right)^2\right)=$$
$$v_{\rho\rho}\sum_{i=1}^Ng_N\left(v_{\rho}^2\left(\frac{x_i}{\rho}\right)^2\right)\left(\frac{x_i}{\rho}\right)^2+\frac{v_{\rho}}{\rho}(N-1)\geq v_{\rho\rho}g_N\left(\frac{v_{\rho}^2}{N}\right)+\frac{v_{\rho}}{\rho}(N-1).$$
The function $v$ solves (\ref{eqnrad}), indeed, $v$ is in $C^2(\overline\omega)$ and it is such that $F(v)\geq0$ and $v>0$ in $\omega$, $v(r)=0$, $v(r/2)\leq\epsilon$ and $v_\rho(z)<0$.\\

c) Since $u,v\in W^{1,2}(\omega)\cap C(\overline{\omega})$, $v$ is a weak subsolution and $u$ is a weak solution to $F(u)=0$, and $v_{|\partial \omega}\leq u_{|\partial \omega}$, applying Lemma \ref{comp}, we obtain that $v\leq u$ in $\omega$. From
$$v_\rho(z)=\frac{\partial v}{\partial\nu}(z)<\frac{\partial u}{\partial\nu}(z),$$
it follows that there exists $x^0\in\omega$ such that $v(x^0)>u(x^0)$, a contradiction.

\endproof

From Hopf's Lemma we derive:

\begin{theorem}[Strong Maximum Principle]\label{thsmpr}

Let $\Omega\subset\mathbb R^N$ be a connected, open and bounded set. Let $u\in W^{1,2}(\Omega)\cap C\left(\overline\Omega\right)$ be a weak supersolution to
$$\sum_{i=1}^N g_i(u_{x_i}^2)u_{x_ix_i}=0.$$
In addition to the assumptions (L) and (G) on $g_i$, assume that $G(\xi)\equiv+\infty$.
Then, if $u$ attains its minimum in $\Omega$, it is a constant.

\end{theorem}

\proof

a) Assume $\min_\Omega u=0$ and set $\mathcal{C}=\{x\in\Omega:u(x)=0\}$. By contradiction, suppose that the open set $\Omega\setminus\mathcal{C}\neq\emptyset$.\\

b) Since $\Omega$ is a connected set, there exist $s\in\mathcal{C}$ and $R>0$ such that $B(s,R)\subset\Omega$ and $B(s,R)\cap(\Omega\setminus\mathcal{C})\neq\emptyset$. Let $p\in B(s,R)\cap(\Omega\setminus\mathcal{C})$. Consider the line $\overline{ps}$. Moving $p$ along this line, we can assume that $B(p,d(p,\mathcal{C}))\subset(\Omega\setminus\mathcal{C})$ and that there exists one point $z\in\mathcal{C}$ such that $d(p,\mathcal{C})=d(p,z)$. Set $r=d(p,\mathcal{C})$. W.l.o.g. suppose that $p=O$.\\

c) The set $\Omega\setminus\mathcal{C}$ satisfies the interior ball condition at $z$, hence Hopf's Lemma implies
$$\frac{\partial u}{\partial\nu}(z)<0.$$
But this is a contradiction: since $u$ attains minimum at $z\in\Omega$, we have that $Du(z)=0$.\\

\endproof

\section{A necessary condition for the validity of Hopf's Lemma}

In this and the following section we consider the operator
\begin{equation}\label{gN}
F(u)=\sum_{i=1}^{N-1}u_{x_ix_i}+g(u_{x_N}^2)u_{x_Nx_N},
\end{equation}
We wish to provide a necessary condition for the validity of Hopf's Lemma in a class of non elliptic equations.

Consider the case
$$G(\xi)=\int^\xi_0{\frac{g(\zeta^2/N)}{\zeta}d\zeta}<+\infty.$$

\begin{theorem}\label{hopfn}

Consider the operator (\ref{gN}), where $g$ satisfies assumptions (L) and (G), and on $(0,\overline t\,]$ (\,$\overline{t}$ defined in assumptions (L)),
$$g^\prime(t)>0 \quad\mbox{and}\quad g(t)+g^\prime(t)t\leq1.$$
If
$$\lim_{t\rightarrow0^+}\frac{\left(g(t)\right)^{3/2}}{tg^\prime(t)}=0,$$
then there exist: an open regular region $\Omega\subset\mathbb R^N$; a radial function
$u\in C^2(\Omega)$ such that $F(u)\leq0$ in $\Omega$ and a point $z\in\partial\Omega$ such that
$u(z)=0$, $u(z)\leq u(x)$ for every $x\in\Omega$ and
$$\frac{\partial u}{\partial\nu}(z)=0$$
where $\nu$ is the outer unit normal to $\Omega$ at $z$.\\

If, in addition, we assume that
\begin{equation}\label{C2}
\frac{g(t)}{tg^\prime(t)} \quad \mbox{is bounded in} \quad (0,\overline t\,]
\end{equation}
then $\Omega$ satisfies the interior ball condition at $z$.

\end{theorem}

\begin{remark}

When $\lim_{t\rightarrow0^+}{g^\prime(t)t}$ exists, it follows that
$$\lim_{t\rightarrow0^+}\frac{\left(g(t)\right)^{3/2}}{tg^\prime(t)}$$
exists, and that
$$g(t)+g^\prime(t)t\leq1\,,\mbox{ on }(0,\overline t\,].$$

\end{remark}
Indeed, we have that
$$\lim_{t\rightarrow0^+}{\left(g(t)+g^\prime(t)t\right)}=0.$$
Otherwise, there exists $K>0$ such that, when $0<t\leq\overline t$,
$g^\prime(t)t\geq K$, so that
$$g(t)t=\int^t_0{\left(g(s)+g^\prime(s)s\right)ds}\geq Kt,$$
and $g(t)\geq K$. From
$$\int^{\xi^2}_0{\frac{g(t)}{t}dt}<+\infty,$$
it follows that $\lim_{t\rightarrow0^+}{g(t)}=0$, a contradiction.\\

The map
$$g(t)=\frac{1}{|\ln(t)|^k},$$
with $k>2$, for $0\leq t\leq1/e$, satisfies the assumption
$$\lim_{t\rightarrow0^+}\frac{\left(g(t)\right)^{3/2}}{tg^\prime(t)}=0.$$\\

The following lemma is instrumental to the proofs of the main results.

\begin{lemma}\label{l1}

Let $g$ satisfies assumptions (L) and (G). Suppose that for every $0<t\leq\overline t$,
$$g^\prime(t)\geq0\quad\mbox{ and }\quad g(t)+g^\prime(t)t\leq1.$$
Set
$$k_1(a)=(1-a)+ag(ta)\quad\mbox{ and }\quad k_2(a)=-a-(1-a)g(ta)$$
For every $0<t\leq\overline t$ (\,$\overline t$ defined in assumptions (L)), $k_1$ and $k_2$ are non increasing in $[0,1]$.

\end{lemma}

\proof

Since, for every $0<t\leq\overline t$,
$$g^\prime(t)\geq0\quad\mbox{ and }\quad g(t)+g^\prime(t)t\leq1,$$
we have that, for every $0\leq a\leq1$,
$$k_1^{\prime}(a)=-1+g(ta)+g^\prime(ta)ta\leq0$$
and
$$k_2^{\prime}(a)=-1+g(ta)-(1-a)g^\prime(ta)t=-1+g(ta)+g^\prime(ta)ta-g^\prime(ta)t\leq0.$$

\endproof

\proof[Proof of Theorem \ref{hopfn}]

a) Let $v$ be a radial function.
Setting $a=\sin^2\theta_{N-1}$, (\ref{gN}) reduces to
$$F(v)=v_{\rho\rho}\left(1-a+ag(v_{\rho}^2a)\right)+\frac{v_{\rho}}{\rho}\left(N-2+a+(1-a)g(v_{\rho}^2a)\right).$$
Let $a=1$, we seek a solution to
\begin{equation}\label{supersol}
v_{\rho\rho}g(v_\rho^2)+(N-1)\frac{v_{\rho}}{\rho}=0
\end{equation}
such that $v_\rho(R(1)+1)=0$ and $v_\rho(\rho)<0$, for every $\rho\in[2,R(1)+1)$.
Consider the Cauchy problem
\begin{equation}\label{csup}
\left\{\begin{array}{l}
        \zeta^\prime=-\dfrac{N-1}{\rho}\dfrac{\zeta}{g(\zeta^2)}\\
        \\
        \zeta(2)=-1.\\
\end{array}
\right.
\end{equation}

We are interested in a negative solution $\zeta$. Define $R(1)$ to be the unique positive real solution to
$$G(-1)-(N-1)\ln\left(\frac{R(1)+1}{2}\right)=0,$$
i.e.
$$R(1)=2e^{\frac{G(-1)}{N-1}}-1.$$
The solution $\zeta$ of (\ref{csup}), satisfies
$$G(\zeta(\rho))-G(-1)=\int^{\zeta(\rho)}_{\zeta(2)}{\frac{g(t^2)}{t}dt}=\int^{\rho}_{2}{-\frac{N-1}{s}ds}=-(N-1)\ln\left(\frac{\rho}{2}\right).$$
Then, for every $\rho\in(2,R(1)+1)$, $G(\zeta(\rho))>0$ and $\zeta(\rho)<0$, while $\zeta(R(1)+1)=0$. Setting $v_{\rho}(\rho)=\zeta(\rho)$ and
$$v(\rho)=\int_{R(1)+1}^{\rho}{v_{\rho}(s)ds},$$
we obtain that $v$ solves (\ref{supersol}) and, for every $\rho\in(2,R(1)+1)$,
$$v_{\rho}(\rho)<v_{\rho}(R(1)+1)=0\quad\mbox{ and }\quad v(\rho)>v(R(1)+1)=0.$$

b) Set, for $\rho\in(1,R(1)]$, $u(\rho)=v(\rho+1).$
Since, for the function $v$, we have
$$v_{\rho\rho}(\rho)g(v_\rho^2(\rho))+(N-1)\frac{v_{\rho}(\rho)}{\rho}=0,$$
at $\rho+1$ we obtain
\begin{equation}\label{B}
u_{\rho\rho}(\rho)g(u_\rho^2(\rho))+(N-1)\frac{u_{\rho}(\rho)}{\rho+1}=0.
\end{equation}
This equality yields, for $\rho\in(1,R(1))$,
$$u_{\rho\rho}g(u_\rho^2)+(N-1)\frac{u_{\rho}}{\rho}=-(N-1)u_{\rho}\left(\frac{1}{\rho+1}-\frac{1}{\rho}\right)<0.$$

c) Let $1/2<a<1$. We wish to find $R(a)\leq R(1)$ such that $u$ is a solution to
$$F(u)=u_{\rho\rho}\left(1-a+ag(u_{\rho}^2a)\right)+\frac{u_{\rho}}{\rho}\left(N-2+a+(1-a)g(u_{\rho}^2a)\right)\leq0,$$
for $\rho\in(1,R(a))$.
Since
$$F(u)=\frac{-u_{\rho}}{\rho(\rho+1)g((u_\rho)^2)}\left[\rho(N-1)\left(1-a+ag((u_{\rho})^2a)\right)-\right.$$
$$\left.(\rho+1)g((u_{\rho})^2)\left(N-2+a+(1-a)g((u_{\rho})^2a)\right)\right],$$
setting
$$k(\rho)=$$
$$\rho(N-1)\left(1-a+ag((u_{\rho})^2a)\right)-(\rho+1)g((u_{\rho})^2)\left(N-2+a+(1-a)g((u_{\rho})^2a)\right),$$
we have that $F(u)\leq0$ if and only if $k(\rho)\leq0$.
From
$$\frac{d}{d\rho}g((u_{\rho})^2)\leq\frac{d}{d\rho}g((u_{\rho})^2a)\leq0,$$
and $N-2+a\geq(N-1)a$ for $N\geq2$, applying Lemma \ref{l1}, we have that
$$k^\prime(\rho)=(N-1)(1-a+ag((u_{\rho})^2a)-g((u_{\rho})^2)(N-2+a+(1-a)g((u_{\rho})^2a))+$$
$$\rho(N-1)a\frac{d}{d\rho}g((u_{\rho})^2a)-(\rho+1)\frac{d}{d\rho}g((u_{\rho})^2)(N-2+a+(1-a)g((u_{\rho})^2a))-$$
$$(\rho+1)g((u_{\rho})^2)(1-a)\frac{d}{d\rho}g((u_{\rho})^2a)\geq$$
$$\rho\left[(N-1)a\frac{d}{d\rho}g((u_{\rho})^2a)-\frac{d}{d\rho}g((u_{\rho})^2)(N-2+a+(1-a)g((u_{\rho})^2a))-\right.$$
$$\left.g((u_{\rho})^2)(1-a)\frac{d}{d\rho}g((u_{\rho})^2a)\right]\geq$$
$$\rho(N-1)a\left(\frac{d}{d\rho}g((u_{\rho})^2a)-\frac{d}{d\rho}g((u_{\rho})^2)\right)\geq0.$$
Since the function $k(\rho)$ is non decreasing, it follows that $F(u)\leq0$, for every $\rho\in(1,R(a))$, if and only if
$$k(R(a))\leq0.$$
We have that
$$k(R(a))=R(a)(N-1)\left(1-a+ag((u_{\rho}(R(a)))^2a)\right)-$$
$$(R(a)+1)g((u_{\rho}(R(a)))^2)\left(N-2+a+(1-a)g((u_{\rho}(R(a)))^2a)\right)\leq$$
$$(N-1)\left[R(a)(1-a)-g((u_{\rho}(R(a)))^2)a\right]\leq$$
$$(N-1)\left[R(1)(1-a)-g((u_{\rho}(R(a)))^2)a\right]\leq$$
$$(N-1)\left[R(1)(1-a)-\frac{g((u_{\rho}(R(a)))^2)}{2}\right].$$
We define $R(a)$ to be a solution to
\begin{equation}\label{A}
R(1)(1-a)-\frac{g((u_{\rho}(R(a)))^2)}{2}=0.
\end{equation}\\

d) In order to solve (\ref{A}) for the unknown $R(a)$, recalling that $1-a=\cos^2\theta_{N-1}=c^2$, let
$$h(r)=\sqrt{\frac{g((u_{\rho}(r))^2)}{2R(1)}}.$$
The function $h$ is decreasing, differentiable and with inverse differentiable.
We have that $|c|=h(R(1-c^2))$, so that $R(1-c^2)=h^{-1}(|c|)$, $R(1-c^2)$ is increasing in $|c|$ and
$$\lim_{c\rightarrow0}R(1-c^2)=R(1).$$
Let $0<|\bar{c}|<1/2$ be such that $R(1-\bar{c}^2)\geq1$, so that, for $c^2\leq\bar{c}^2$, we have $R(1-c^2)\geq R(1-\bar{c}^2)\geq1$. We have obtained that, for every $1\geq a\geq 1-|c|^2$, there exists $R(a)$ such that (\ref{A}) holds.
It follows that
$$k(R(a))\leq\left[R(1)(1-a)-\frac{g((u_{\rho}(R(a)))^2)}{2}\right]=0,$$
so that the function $u$ solves $F(u)\leq0$ for every $\rho\in(1,R(a))$.\\

e) Set $\Omega=\{(x_1,\dots,x_N)\in\mathbb R^N:\rho\in(1,R(1-c^2))\mbox{ and }|c|<|\bar c|\}$. $\Omega\subset\mathbb R^N$ is a connected, open and bounded set and $u\in W^{1,2}(\Omega)\cap C\left(\overline\Omega\right)$ is a weak solution to
$F(u)\leq0.$ The point $z=(R,0,\dots,0)\in\partial\Omega$ is such that $u(z)<u(x)$, for all $x\in\Omega$. We wish to show that $\Omega$ is regular in a neighborhood of $z=(R,0,\dots,0)$.
Since $\frac{d}{dc}R(1-c^2)$ exists, in $(0,|\bar c|)$, to prove our claim it is sufficient to show that
$$\lim_{c\rightarrow0}\frac{d}{dc}R(1-c^2)=0.$$
Recalling (\ref{B}), we have that
$$\frac{d}{dc}R(1-c^2)=\left(h^{-1}(c)\right)^\prime=\frac{1}{h^\prime(R(1-c^2))}=$$
\begin{equation}\label{dR}
-\frac{\sqrt{2}R(1)}{N-1}\left(R(1-c^2)+1\right)\frac{\left(g((u_{\rho}(R(1-c^2)))^2)\right)^{3/2}}{(u_{\rho}(R(1-c^2)))^2g^\prime((u_{\rho}(R(1-c^2)))^2)}.
\end{equation}
Since
$$\lim_{c\rightarrow0}u_{\rho}(R(1-c^2))=0$$
and
$$\lim_{t\rightarrow0^+}\frac{\left(g(t)\right)^{3/2}}{tg^\prime(t)}=0,$$
it follows that
$$\lim_{c\rightarrow0}\frac{d}{dc}R(1-c^2)=0.$$
Since ${\frac{d}{d\theta}\cos\theta_{N-1}}_{|\theta_{N-1}=\frac{\pi}{2}}=1$, this shows that
$\frac{d}{d\theta}R(1-\cos^2\theta_{N-1})=0$, and $\Omega$ is regular.\\

f) To prove the validity of the interior ball condition at $z=(R,0,\dots,0)$, it is enough to show that the second derivative of $R(1-c^2)$ is bounded at $c=0$, i.e. that
$$\left|\frac{1}{c}\frac{d}{dc}R(1-c^2)\right|$$
is bounded. Set
$$t(c)=(u_{\rho}(R(1-c^2)))^2,$$
from (\ref{dR}) we obtain
$$\frac{d}{dc}R(1-c^2)=-\frac{\sqrt{2}R(1)}{N-1}\left(R(1-c^2)+1\right)\frac{(g(t(c)))^{3/2}}{t(c)g^\prime(t(c))},$$
and from (\ref{B})
$$\frac{dt(c)}{dc}=2u_{\rho}(R(1-c^2))u_{\rho\rho}(R(1-c^2))\frac{d}{dc}(R(1-c^2))=$$
$$2\sqrt{2}R(1)\frac{(g(t(c)))^{1/2}}{g^\prime(t(c))}$$
and
$$\frac{d}{dc}(g(t(c)))^{1/2}=\frac{g^\prime(t(c))}{2(g(t(c)))^{1/2}}\frac{dt(c)}{dc}=\sqrt{2}R(1).$$
From $g(t(0))=0$, we obtain that
$$(g(t(c)))^{1/2}=\sqrt{2}R(1)c$$
and
$$\frac{(g(t(c)))^{3/2}}{ct(c)g^\prime(t(c))}=\sqrt{2}R(1)\frac{g(t(c))}{t(c)g^\prime(t(c))}.$$
From condition (\ref{C2}) we obtain
$$\left|\frac{1}{c}\frac{d}{dc}R(1-c^2)\right|\leq M.$$

\begin{figure}
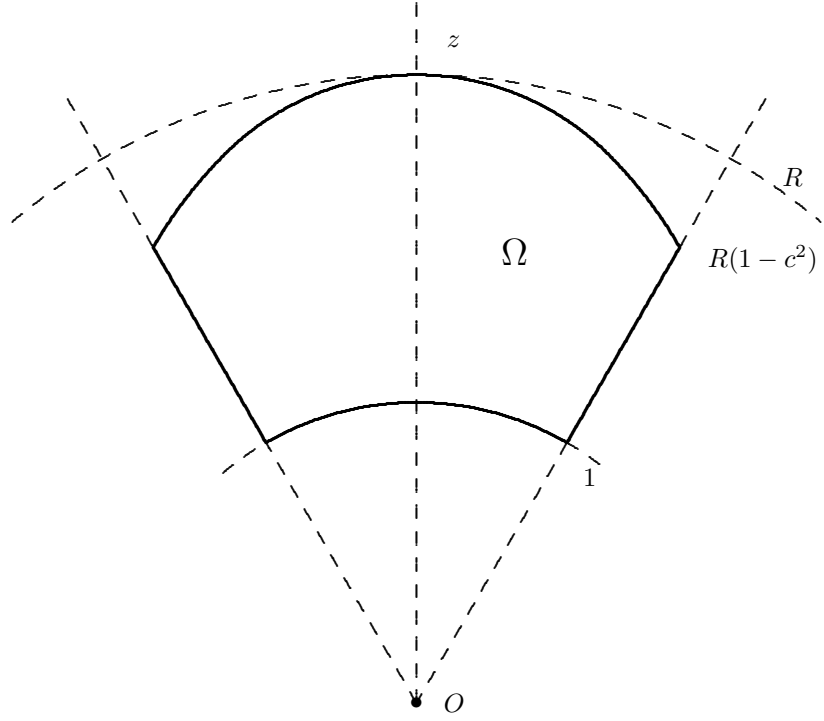

\begin{center}
\beginpicture

\setcoordinatesystem units < 1cm, 1cm>
\setplotarea x from -7 to  7, y from -1 to  9
\put {$\bullet$} at  0  0
\setplotsymbol({\large.})
\plot   2  3.464102     3.5  6.062178 /
\plot  -2  3.464102    -3.5  6.062178 /
\circulararc -60 degrees from -2  3.464102 center at  0  0
\circulararc -15 degrees from  2.546488  7.304819 center at -1.696153  3.062178
\circulararc -14.99999 degrees from -3.5  6.062178 center at  1.696151  3.062177
\circulararc -89.99998 degrees from -2.546491  7.304819 center at -0.00000143  4.75833
\setplotsymbol({\tiny.})
\setdashes
\plot   0  0    -0.00000143  9.35961 /
\plot   2  3.464102     0  0 /
\plot   0  0    -2  3.464102 /
\plot   3.5  6.062178     4.679805  8.105659 /
\plot  -3.5  6.062178    -4.679805  8.105659 /
\circulararc -80 degrees from -5.373454  6.403832 center at  0  0
\circulararc -10.00001 degrees from  2  3.464102 center at  0  0
\circulararc -9.99999 degrees from -2.571151  3.064178 center at  0  0
\setsolid
\put {$O$} at  .5  0
\put {$z$} at  .49  8.8
\put {$1$} at  2.3  3
\put {$R$} at  5  7
\put {$R(1-c^2)$} at  4.6  5.9
\put {\LARGE{$\Omega$}} at  1.3  6

\endpicture
\caption{$\Omega$ in the case $N=2$.}
\label{fignec}
\end{center}
\end{figure}

\endproof

\section{A sufficient condition for the validity of the Strong Maximum Principle}

Consider the case
$$G(\xi)=\int^\xi_0{\frac{g(\zeta^2/N)}{\zeta}d\zeta}<+\infty.$$

We wish to prove the following theorem.

\begin{theorem}[The Strong Maximum Principle]\label{thsmp}

Let $\Omega\subset\mathbb R^N$ be a connected, open and bounded set. Let $u\in W^{1,2}(\Omega)\cap C\left(\overline\Omega\right)$ be a weak supersolution to
$$\sum_{i=1}^{N-1}u_{x_ix_i}+g(u_{x_N}^2)u_{x_Nx_N}=0,$$
where $g$ satisfies assumptions (L) and (G) and, on $(0,\overline t\,]$ (\,$\overline{t}$ defined in assumptions (L)),
$$g^\prime(t)\geq0\quad\mbox{ and }\quad g(t)+g^\prime(t)t\leq1.$$
Moreover, suppose that there exists $K>0$ such that, for every $0<\xi^2/N\leq\overline{t}$ , we have
$$\sqrt{g(\xi^2/N)}\leq K\left(e^{G(\xi)}-1\right).$$
Then, if $u$ attains its minimum in $\Omega$, it is a constant.

\end{theorem}

\begin{remark}

When the function $g$ satisfies the condition
$$g^\prime(t)t\leq2K\left(g(t)\right)^{3/2},$$
for every $0<t\leq\overline t$, then it satisfies
$$g(t)+g^\prime(t)t\leq1$$
and
$$\sqrt{g(\xi^2/N)}\leq K\left(e^{G(\xi)}-1\right),$$
for every $0<\xi^2/N\leq\overline t$.

\end{remark}
Indeed,
since $G(\xi)<+\infty$, we have that $\lim_{t\rightarrow0}g(t)=0$. Hence, we can suppose that $g(t)\leq1/(2K+1)$, for $0<t\leq\overline t$, so that
$$g(t)+g^\prime(t)t\leq g(t)+2K\left(g(t)\right)^{3/2}\leq1.$$
Moreover, since $g(0)=0$, $G(0)=0$ and
$$\left(\sqrt{g(\xi^2/N)}\right)^\prime\leq2K\frac{g(\xi^2/N)}{\xi}\leq  K\left(e^{G(\xi)}-1\right)^\prime,$$
we obtain that
$$\sqrt{g(\xi^2/N)}\leq K\left(e^{G(\xi)}-1\right).$$

\begin{remark}

Among the functions $g$ such that
$$\lim_{t\rightarrow0^+}\frac{\left(g(t)\right)^{3/2}}{tg^\prime(t)}$$
exists, there exists $K>0$ such that, for every $0<\xi^2/N\leq\overline{t}$, we have
$$\sqrt{g(\xi^2/N)}\leq K\left(e^{G(\xi)}-1\right)$$
if and only if
$$\lim_{t\rightarrow0^+}\frac{\left(g(t)\right)^{3/2}}{tg^\prime(t)}>0.$$

\end{remark}

An example of a map satisfying the assumptions of the theorem above, is given by
$$g(t)=\frac{1}{(\ln(t))^2},$$
for $0\leq t\leq1/e^4$. For $0\leq\xi^2/N\leq1/e^4$, we have
$$G(\xi)=\int^{\xi}_{0}{\frac{1}{\zeta(\ln(\zeta^2/N))^2}d\zeta}=-\frac{1}{2\ln(\xi^2/N)},$$
and
$$\sqrt{g(\xi^2/N)}=\frac{1}{|\ln(\xi^2/N)|}\leq2\left(e^{\frac{1}{2|\ln(\xi^2/N)|}}-1\right)=e^{G(\xi)}-1.$$

Set
$$\mathcal R(\lambda,\lambda_N)=\{(x_1,\dots,x_N):|x_i|\leq\lambda,\mbox{ for }i=1,\dots,N-1, |x_N|\leq\lambda_N\}.$$

To the opposite of the proof of Lemma \ref{hopfr} and Theorem \ref{thsmpr}, we will build a subsolution that {\it is not} radially symmetric. This construction is provided by next theorem.

\begin{theorem}\label{subsol}

Under the same assumptions on $g$ as on Theorem \ref{thsmp}, for every $r>0$ and every $\epsilon$, there exist: $l,l_N$; an open convex region $\mathcal{A}\subset\mathcal R(l,l_N)$; a function $v\in W^{1,2}(\omega)\cap C^1(\omega)\cap C(\overline{\omega})$, where $\omega=B(\mathcal{A},r)\setminus\overline{\mathcal{A}}$, such that

i) $0\leq l\leq2Kr$, and $0\leq l_N\leq r/4$;

ii)\begin{equation}\label{eqnN}
\left\{\begin{array}{ll}
        v \mbox{ is a weak solution to }F(v)\geq0 & \mbox{ in } \omega\\
        v>0 & \mbox{ in }\omega\\
        v=0 & \mbox{ in }\partial B(\mathcal{A},r)\\
        v\leq\epsilon & \mbox{ in }\partial\mathcal{A}.\\
        \end{array}
\right.
\end{equation}

\end{theorem}

\proof

Fix $r$; we can assume that $\epsilon$ is such that $0<\epsilon^2/r^2\leq\overline t$ and that $$2\sqrt{g\left(\frac{\epsilon^2}{r^2}\right)}+\frac{1}{2}\sqrt{g\left(\frac{\epsilon^2}{r^2}\right)}\left|\ln g\left(\frac{\epsilon^2}{r^2}\right)\right|\leq\frac{1}{4K}.$$
Fix the origin $O^0=(0,\dots,0)$, and set polar coordinates as
$$\left\{\begin{array}{l}
        x_1=\rho\cos\theta_{N-1}\dots\cos\theta_2\cos\theta_1\\
        x_2=\rho\cos\theta_{N-1}\dots\cos\theta_2\sin\theta_1\\
        \dots\\
        x_N=\rho\sin\theta_{N-1}.
        \end{array}
\right.$$\\

1) When $w$ is a radial function, setting $a=\sin^2\theta_1$ , $F$ reduces to
$$F(w)=w_{\rho\rho}\left(1-a+ag(w_{\rho}^2a)\right)+\frac{w_{\rho}}{\rho}\left(N-2+a+(1-a)g(w_{\rho}^2a)\right).$$
For $a=1$, we seek a solution to
\begin{equation}\label{pN}
(N-1)\frac{w_{\rho}}{\rho}+w_{\rho\rho}g(w_\rho^2)=0.
\end{equation}
such that $w_\rho(R(1))=-\epsilon/r$ and $w_\rho(\rho)<0$, for every $\rho\in[R(1),R(1)+r)$.
Consider the Cauchy problem
\begin{equation}\label{c2}
\left\{\begin{array}{l}
        \zeta^\prime=-\dfrac{N-1}{\rho}\dfrac{\zeta}{g(\zeta^2)}\\
        \\
        \zeta(R(1))=-\dfrac{\epsilon}{r}.\\
\end{array}\right.
\end{equation}
We are interested in a negative solution $\zeta$.
Define $R(1)$ to be the unique positive real solution to
$$G(-\epsilon/r)-(N-1)\ln\left(\frac{R(1)+r}{R(1)}\right)=0,$$
i.e.
$$R(1)=\frac{r}{e^{\frac{G(-\epsilon/r)}{N-1}}-1}.$$
Consider the unique solution $\zeta$ of (\ref{c2}), such that $\zeta(R(1))=-\epsilon/r$, i.e., such that
$$G(\zeta(\rho))-G(-\epsilon/r)=\int^{\zeta(\rho)}_{\zeta(R(1))}{\frac{g(t^2)}{t}dt}=\int^{\rho}_{R(1)}{-\frac{N-1}{s}ds}=-(N-1)\ln\left(\frac{\rho}{R(1)}\right).$$
Then, for every $\rho\in[R(1),R(1)+r)$, $G(\zeta(\rho))>0$ and $\zeta(\rho)<0$, while $\zeta(R(1)+r)=0$. Setting $w_{\rho}(\rho)=\zeta(\rho)$ and
$$w(\rho)=\int_{R(1)+r}^{\rho} {w_{\rho}(s)ds},$$
we obtain that $w$ solves (\ref{pN}) and, for every $\rho\in(R(1),R(1)+r)$,
$$-\epsilon/r=w_{\rho}(R(1))<w_{\rho}(\rho)<w_{\rho}(R(1)+r)=0$$
and
$$0=w(R(1)+r)<w(\rho)<w(R(1))\leq\epsilon.$$

2) Applying Lemma \ref{l1}, we infer that the function $w$ defined in 1) is actually a solution to
$$F(w)=w_{\rho\rho}(1-a)+\frac{w_{\rho}}{\rho}(N-2+a)+g(w_\rho^2a)\left(w_{\rho\rho}a+\frac{w_{\rho}}{\rho}(1-a)\right)\geq0,$$
for every $0\leq a\leq1$ and every $\rho\in(R(1),R(1)+r)$.\\

3) Let $\bar{a}<1$. We wish to find the smallest $R(\bar{a})>0$ such that, setting
$$w^{\bar{a}}(\rho)=w(\rho-R(\bar{a})+R(1)),$$
the function $w^{\bar{a}}$ is a solution to
\begin{equation}\label{eqa}
F(w^{\bar{a}})=w^{\bar{a}}_{\rho\rho}(1-\bar{a})+\frac{w^{\bar{a}}_{\rho}}{\rho}(N-2+\bar{a})+g((w^{\bar{a}}_\rho)^2\bar{a})\left(w^{\bar{a}}_{\rho\rho}\bar{a}+\frac{w^{\bar{a}}_{\rho}}{\rho}(1-\bar{a})\right)\geq0,
\end{equation}
for every $\rho\in(R(\bar{a}),R(\bar{a})+r)$.

Since, for the function $w$, we have
$$(N-1)\frac{w_{\rho}(\rho)}{\rho}+g(w_\rho^2(\rho))w_{\rho\rho}(\rho)=0,$$
at $\rho-R(\bar{a})+R(1)$ we obtain
$$(N-1)\frac{w^{\bar{a}}_{\rho}(\rho)}{\rho-R(\bar{a})+R(1)}+g((w^{\bar{a}}_\rho(\rho))^2) w^{\bar{a}}_{\rho\rho}(\rho)=0.$$
This equality yields
$$(N-2+\bar{a})\frac{w^{\bar{a}}_{\rho}}{\rho}+\bar{a}g((w^{\bar{a}}_{\rho})^2\bar{a})w^{\bar{a}}_{\rho\rho}=$$
$$\frac{w^{\bar{a}}_{\rho}}{\rho(\rho-R(\bar{a})+R(1))}\left((N+2-\bar{a})(\rho-R(\bar{a})+R(1))
-\bar{a}\rho(N-1)\frac{g((w^{\bar{a}}_{\rho})^2\bar{a})}{g((w^{\bar{a}}_{\rho})^2)}\right)$$
and
$$(1-\bar{a})\left(w^{\bar{a}}_{\rho\rho}+\frac{w^{\bar{a}}_{\rho}}{\rho}g((w^{\bar{a}}_{\rho})^2\bar{a})\right)=$$
$$\frac{(1-\bar{a})w^{\bar{a}}_{\rho}}{\rho(\rho-R(\bar{a})+R(1))}\left((\rho-R(\bar{a})+R(1))
g((w^{\bar{a}}_{\rho})^2\bar{a})-\frac{(N-1)\rho}{g((w^{\bar{a}}_{\rho})^2)}\right).$$

Since
$$F(w^{\bar{a}})=(N-2+\bar{a})\frac{w^{\bar{a}}_{\rho}}{\rho}+\bar{a}g((w^{\bar{a}}_{\rho})^2\bar{a})w^{\bar{a}}_{\rho\rho}+(1-\bar{a})\left(w^{\bar{a}}_{\rho\rho}+\frac{w^{\bar{a}}_{\rho}}{\rho}g((w^{\bar{a}}_{\rho})^2)\right),$$
we obtain that
$F(w^{\bar{a}})\geq0$ if and only if
$$(\rho-R(\bar{a})+R(1))\left(N-2+\bar{a}+(1-\bar{a})g((w^{\bar{a}}_{\rho})^2\bar{a})\right)-$$
$$\frac{(N-1)\rho}{g((w^{\bar{a}}_{\rho})^2)}\left(1-\bar{a}+\bar{a}g((w^{\bar{a}}_{\rho})^2\bar{a})\right)\leq0.$$
Set
$$k(\rho)=(\rho-R(\bar{a})+R(1))\left(N-2+\bar{a}+(1-\bar{a})g((w^{\bar{a}}_{\rho})^2\bar{a})\right)-$$
$$\frac{(N-1)\rho}{g((w^{\bar{a}}_{\rho})^2)}\left(1-\bar{a}+\bar{a}g((w^{\bar{a}}_{\rho})^2\bar{a})\right)$$
Since $$\frac{d}{d\rho}g((w^{\bar{a}}_{\rho})^2)\leq\frac{d}{d\rho}g((w^{\bar{a}}_{\rho})^2\bar a)\leq0,$$
applying Lemma \ref{l1}, we have that
$$k^\prime(\rho)=$$
$$\left(N-2+\bar{a}+(1-\bar{a})g((w^{\bar{a}}_{\rho})^2\bar{a})\right)+
(\rho-R(\bar{a})+R(1))(1-\bar{a})\frac{d}{d\rho}g((w^{\bar{a}}_{\rho})^2\bar a)-$$
$$(N-1)\left(\frac{1}{g((w^{\bar{a}}_{\rho})^2)}-\frac{\rho}{(g((w^{\bar{a}}_{\rho})^2))^2}\frac{d}{d\rho}g((w^{\bar{a}}_{\rho})^2)\right)\left(1-\bar{a}+\bar{a}g((w^{\bar{a}}_{\rho})^2\bar{a})\right)-$$
$$(N-1)\frac{\rho}{g((w^{\bar{a}}_{\rho})^2)}\frac{d}{d\rho}g((w^{\bar{a}}_{\rho}(\rho))^2\bar a)\bar{a}\leq$$
$$\frac{(N-1)\rho}{(g((w^{\bar{a}}_{\rho})^2))^2}\frac{d}{d\rho}g((w^{\bar{a}}_{\rho})^2)\left(1-\bar{a}+\bar{a}g((w^{\bar{a}}_{\rho})^2\bar{a})\right)-\frac{(N-1)\rho}{g((w^{\bar{a}}_{\rho})^2)}\frac{d}{d\rho}g((w^{\bar{a}}_{\rho}(\rho))^2\bar a)\bar{a}\leq$$
$$\frac{(N-1)\rho}{(g((w^{\bar{a}}_{\rho})^2))^2}\frac{d}{d\rho}g((w^{\bar{a}}_{\rho})^2)\left(1-\bar{a}+\bar{a}g((w^{\bar{a}}_{\rho})^2\bar{a})-\bar{a}g((w^{\bar{a}}_{\rho})^2)\right)\leq$$
$$\frac{(N-1)\rho}{(g((w^{\bar{a}}_{\rho})^2))^2}(1-\bar{a})g((w^{\bar{a}}_{\rho})^2)\frac{d}{d\rho}g((w^{\bar{a}}_{\rho})^2)\leq0.$$
Since the function $k(\rho)$ is non increasing, it follows that $F(w^{\bar{a}})\geq0$, for every $\rho\in(R(\bar{a}),R(\bar{a})+r)$, if and only if
$$R(1)\left(N-2+\bar{a}+(1-\bar{a})g((w^{\bar{a}}_{\rho}(R(\bar{a})))^2\bar{a})\right)-$$
$$\frac{(N-1)R(\bar{a})}{g((w^{\bar{a}}_{\rho}(R(\bar{a}))^2)}\left(1-\bar{a}+\bar{a}g((w^{\bar{a}}_{\rho}(R(\bar{a})))^2\bar{a})\right)=$$
$$R(1)\left(N-2+\bar{a}+(1-\bar{a})g\left(\frac{\epsilon^2}{r^2}\bar{a}\right)\right)-
\frac{(N-1)R(\bar{a})}{g\left(\frac{\epsilon^2}{r^2}\right)}\left(1-\bar{a}+\bar{a}g\left(\frac{\epsilon^2}{r^2}\bar{a}\right)\right)\leq0,$$
i.e. if and only if
$$R(\bar{a})\geq g\left(\frac{\epsilon^2}{r^2}\right)\frac{R(1)}{N-1}\frac{N-2+\bar{a}+g\left(\frac{\epsilon^2}{r^2}\bar{a}\right)(1-\bar{a})}{1-\bar{a}+g\left(\frac{\epsilon^2}{r^2}\bar{a}\right)\bar{a}}.$$
Hence, we define
$$R(\bar{a})=g\left(\frac{\epsilon^2}{r^2}\right)\frac{R(1)}{N-1}\frac{N-2+\bar{a}+g\left(\frac{\epsilon^2}{r^2}\bar{a}\right)(1-\bar{a})}{1-\bar{a}+g\left(\frac{\epsilon^2}{r^2}\bar{a}\right)\bar{a}}.$$\\

4) The function $w^{\bar{a}}$ defined in point 3) is a solution to
$$F(w^{\bar{a}})=w^{\bar{a}}_{\rho\rho}(1-a)+\frac{w^{\bar{a}}_{\rho}}{\rho}(N-2+a)+
g((w^{\bar{a}}_\rho)^2a)\left(w^{\bar{a}}_{\rho\rho}a+\frac{w^{\bar{a}}_{\rho}}{\rho}(1-a)\right)\geq0,$$
for every $a<\bar{a}$. Indeed, applying Lemma \ref{l1} we obtain that, for every $\rho\in(R(\bar{a}),R(\bar{a})+r)$,
$$(\rho-R(\bar{a})+R(1))\left(N-2+a+(1-a)g((w^{\bar{a}}_{\rho})^2a)\right)-$$
$$\frac{(N-1)\rho}{g((w^{\bar{a}}_{\rho})^2)}\left(1-a+ag((w^{\bar{a}}_{\rho})^2a)\right)\leq$$
$$(\rho-R(\bar{a})+R(1))\left(N-2+\bar{a}+(1-\bar{a})g((w^{\bar{a}}_{\rho})^2\bar{a})\right)-$$
$$\frac{(N-1)\rho}{g((w^{\bar{a}}_{\rho})^2)}\left(1-\bar{a}+\bar{a}g((w^{\bar{a}}_{\rho})^2\bar{a})\right)\leq0.$$ \\

5) Assume we have a partition $\alpha$ of $[0,\pi/2]$, $\alpha=\{0=\alpha_n<\dots<\alpha_1<\alpha_0=\pi/2\}$. This partition defines two partitions of $[0,1]$, given by $c_i=\cos\alpha_i$ and $s_i=\sin\alpha_i$.

Consider the sums
$$S_1(\alpha)=\sum_{i=1}^{n}\left(R(1-c_{i-1}^2)-R(1-c_i^2)\right)c_i=R(1)c_1+\sum_{i=1}^{n-1}R(1-c_i^2)\left(c_{i+1}-c_i\right)$$
and
$$S_2(\alpha)=\sum_{i=1}^{n}R(s_{i-1}^2)(s_{i-1}-s_i)=R(1)(1-s_1)+\sum_{i=1}^{n-1}R(s_i^2)(s_i-s_{i+1}),$$
where, in the previous equalities, we have taken into account that $R(1-c_n^2)=R(0)=0$. Our purpose is to provide a partition $\alpha$ and corresponding estimates for $S_1(\alpha)$ and $S_2(\alpha)$ that are independent of $\epsilon$.

The sums
$$\sum_{i=1}^{n-1}R(1-c^2_i)(c_{i+1}-c_i)\quad\mbox{ and }\quad\sum_{i=1}^{n-1}R(s_i^2)(s_i-s_{i+1})$$
are Riemann sums for the integrals
$$\int^1_{c_1}R(1-c^2)dc\quad\mbox{ and }\quad\int^{s_1}_0R(s^2)ds.$$
Consider the first integral. From
$$\frac{N-1-c^2+g\left(\frac{\epsilon^2}{r^2}(1-c^2)\right)c^2}{c^2+g\left(\frac{\epsilon^2}{r^2}(1-c^2)\right)(1-c^2)}\leq\frac{N-1}{c^2}$$
we obtain that
$$\int^1_{c_1}R(1-c^2)dc=g\left(\frac{\epsilon^2}{r^2}\right)\frac{R(1)}{N-1}\int^1_{c_1}\frac{N-1-c^2+g\left(\frac{\epsilon^2}{r^2}(1-c^2)\right)c^2}{c^2+g\left(\frac{\epsilon^2}{r^2}(1-c^2)\right)(1-c^2)}dc\leq$$
$$R(1)g\left(\frac{\epsilon^2}{r^2}\right)\int^1_{c_1}\frac{dc}{c^2}.$$
Set
$$S_x^*(c)=R(1)c+R(1)g\left(\frac{\epsilon^2}{r^2}\right)\int^1_c{\frac{db}{b^2}}=R(1)\left(c+g\left(\frac{\epsilon^2}{r^2}\right)\left(\frac{1}{c}-1\right)\right)=$$
$$R(1)g\left(\frac{\epsilon^2}{r^2}\right)\left(\frac{c}{g\left(\frac{\epsilon^2}{r^2}\right)}+\frac{1}{c}-1\right).$$
Evaluating the last term at the minimum point $c=\sqrt{g\left(\frac{\epsilon^2}{r^2}\right)},$ we obtain
$$S_x^*\left(\sqrt{g\left(\frac{\epsilon^2}{r^2}\right)}\right)=R(1)g\left(\frac{\epsilon^2}{r^2}\right)\left(\frac{2}{\sqrt{g\left(\frac{\epsilon^2}{r^2}\right)}}-1\right)=$$
$$\frac{2r\sqrt{g\left(\frac{\epsilon^2}{r^2}\right)}}{e^{G(\epsilon/r)}-1}-R(1)g\left(\frac{\epsilon^2}{r^2}\right).$$
We fix $c_1=\sqrt{g\left(\frac{\epsilon^2}{r^2}\right)}$, so that $\alpha_1=\arccos{\sqrt{g\left(\frac{\epsilon^2}{r^2}\right)}}$.

Consider the second integral
$$\int^{s_1}_0R(s^2)ds.$$
From
$$\frac{N-2+s^2+g\left(\frac{\epsilon^2}{r^2}s^2\right)(1-s^2)}{1-s^2+g\left(\frac{\epsilon^2}{r^2}\right)s^2}\leq\frac{N-2+s^2+g\left(\frac{\epsilon^2}{r^2}s^2\right)(1-s^2)}{1-s^2}$$
we obtain that
$$\int^{s_1}_0R(s^2)ds=g\left(\frac{\epsilon^2}{r^2}\right)\frac{R(1)}{N-1}\int^{s_1}_0
\frac{N-2+s^2+g\left(\frac{\epsilon^2}{r^2}s^2\right)(1-s^2)}{1-s^2+g\left(\frac{\epsilon^2}{r^2}s^2\right)s^2}ds\leq$$
$$g\left(\frac{\epsilon^2}{r^2}\right)\frac{R(1)}{N-1}\int^{s_1}_0\frac{N-2+s^2+g\left(\frac{\epsilon^2}{r^2}\right)(1-s^2)}{1-s^2}ds.$$
Set
$$S_y^*(s)=R(1)(1-s)+g\left(\frac{\epsilon^2}{r^2}\right)\frac{R(1)}{N-1}\int^s_0\frac{N-2+b^2+g\left(\frac{\epsilon^2}{r^2}\right)(1-b^2)}{1-b^2}db=$$
$$R(1)(1-s)+g\left(\frac{\epsilon^2}{r^2}\right)\frac{R(1)}{N-1}\left[\left(g\left(\frac{\epsilon^2}{r^2}\right)-1\right)s+\frac{N-1}{2}\ln\left(\frac{1+s}{1-s}\right)\right].$$
Since
$$1-\sqrt{1-g\left(\frac{\epsilon^2}{r^2}\right)}\leq g\left(\frac{\epsilon^2}{r^2}\right),$$
evaluating the last term at the point $s_1=\sin\alpha_1=\sqrt{1-g\left(\frac{\epsilon^2}{r^2}\right)}$, we obtain
$$S_x^*\left(\sqrt{1-g\left(\frac{\epsilon^2}{r^2}\right)}\right)=R(1)\left(1-\sqrt{1-g\left(\frac{\epsilon^2}{r^2}\right)}\right)+$$
$$g\left(\frac{\epsilon^2}{r^2}\right)\frac{R(1)}{N-1}\left[\left(g\left(\frac{\epsilon^2}{r^2}\right)-1\right)\sqrt{1-g\left(\frac{\epsilon^2}{r^2}\right)}-\frac{1}{2}\ln\left(\frac{1-\sqrt{1-g\left(\frac{\epsilon^2}{r^2}\right)}}{1+\sqrt{1-g\left(\frac{\epsilon^2}{r^2}\right)}}\right)\right]\leq$$
$$R(1)\sqrt{g\left(\frac{\epsilon^2}{r^2}\right)}\left(\sqrt{g\left(\frac{\epsilon^2}{r^2}\right)}+\frac{1}{2}\sqrt{g\left(\frac{\epsilon^2}{r^2}\right)}\left|\ln g\left(\frac{\epsilon^2}{r^2}\right)\right|\right).$$
To define the other points of the required partition $\alpha$, consider the integrals
$$\int^1_{c_1}{R(1-c^2)dc}\quad\mbox{ and }\quad\int^{s_1}_0R(s^2)ds.$$
Set
$$\sigma=R(1)g\left(\frac{\epsilon^2}{r^2}\right).$$
By the basic theorem of Riemann integration, taking a partition $\alpha$ with mesh size small enough, the value of the Riemann sums
$$\sum_{i=1}^{n-1}R(1-c^2_i)(c_{i+1}-c_i)\quad\mbox{ and }\quad\sum_{i=1}^{n-1}R(s^2_i)(s_i-s_{i+1})$$
differs from
$$\int^1_{c_1}{R(1-c^2)dc}\quad\mbox{ and }\quad\int^{s_1}_0R(s^2)ds$$
by less than $\sigma$. In particular we obtain
$$S_1(\alpha)=R(1)c_1+\sum_{i=1}^{n-1}R(1-c^2_i)\left(c_{i+1}-c_i\right)\leq R(1)c_1+\int^1_{c_1}{R(1-c^2)dc}+\sigma\leq$$
$$\frac{2r\sqrt{g\left(\frac{\epsilon^2}{r^2}\right)}}{e^{\frac{G(\epsilon/r)}{N-1}}-1}\leq2Kr$$
and
$$S_2(\alpha)=R(1)(1-s_1)+\sum_{i=1}^{n-1}R(s^2_i)(s_i-s_{i+1})\leq$$
$$R(1)(1-s_1)+\int^{s_1}_0R(s^2)ds+\sigma\leq$$
$$\frac{r\sqrt{g\left(\frac{\epsilon^2}{r^2}\right)}}{e^{\frac{G(\epsilon/r)}{N-1}}-1}\left(2\sqrt{g\left(\frac{\epsilon^2}{r^2}\right)}+\frac{1}{2}\sqrt{g\left(\frac{\epsilon^2}{r^2}\right)}\left|\ln g\left(\frac{\epsilon^2}{r^2}\right)\right|\right)\leq\frac{r}{4}.$$\\

6) With respect to the coordinates fixed at the beginning of the proof, consider $x_i\geq0$. Set
$$\mathcal{D}_0=\{(x_1,\dots,x_N):R(1)<\rho<R(1)+r,\sqrt{a_1}\leq\sin\theta_{N-1}\leq\sqrt{a_0}=1\}$$
and on $\mathcal{D}_0$ define the function
$$v_0(x_1,\dots,x_N)=w(\rho(x_1,\dots,x_N)).$$
By point 1), the function $v_0$ is of class $C^2(int(\mathcal{D}_0))$ and satisfies, pointwise, the inequality $F(v_0)\geq0$.
Having defined $v_0$, define $v_1$ as follows. Set:
$$O_1(\theta_1,\dots,\theta_{N-2})=(O^1_{x_1},\dots,O^1_{x_N})=$$
$$(R(1)-R(a_1))\left(\sqrt{1-a_1}\cos\theta_{N-2}\dots\cos\theta_1,\sqrt{1-a_1}\cos\theta_{N-2}\dots\sin\theta_1,\dots,\sqrt{a_1}\right),$$
$$\rho_1(x_1,\dots,x_N)=\sqrt{(x_1-O^1_{x_1})^2+\dots+(x_N-O^1_{x_N})^2}$$
and
$$\sin\theta_{N-1}^1(x_1,\dots,x_N)=\frac{x_N-O^1_{x_N}}{\rho_1(x_1,\dots,x_N)}.$$
Recalling the definition of $w^{a_1}$ in 3), consider
$$\mathcal{D}_1=\left\{(x_1,\dots,x_N):R(a_1)<\rho_1(x_1,\dots,x_N)<R(a_1)+r,\right.$$
$$\left.\sqrt{a_2}\leq\sin\theta_{N-1}^1(x_1,\dots,x_N)\leq\sqrt{a_1}\right\}$$
and, on $\mathcal{D}_1$, set
$$v_1(x_1,\dots,x_N)=w^{a_1}(\rho_1(x_1,\dots,x_N)).$$
The function $v_1$ is of class $C^2(int(\mathcal{D}_1))$. We claim that $v_1$ still satisfies $F(v_1)\geq0$.
Remark that the set of the points $O^1$ is equal to
$$\mathcal O^1=\left\{(x_1,\dots,x_N):\rho=R(1)-R(a_1),\quad\sin\theta_{N-1}=\sqrt{a_1}\right\}$$
and that for every point $p\in\mathcal{D}_1$, the corresponding point $O^1(p)$ is the projection of $p$ on $\mathcal O^1$, while $\rho_1(p)=d(p,\mathcal O^1)$. Then we obtain
$$\frac{\partial\rho_1}{\partial\theta_i}=0 \quad\mbox{for every}\quad i=1,\dots,N-2,$$
$$\frac{\partial O^1_{x_i}}{\partial x_i}\geq0 \quad\mbox{for every}\quad i=1,\dots,N-1,$$
$$\frac{\partial O^1_{x_N}}{\partial x_N}=0.$$
and
$$\nabla v_1=\frac{w^{a_1}_{\rho}(\rho_1)}{\rho_1} \left(x_1-O^1_{x_1},\dots,x_N-O^1_{x_N}\right)=$$
$$w^{a_1}_{\rho}(\rho_1)\left(\cos\theta_{N-1}^1\dots\cos\theta_1,\cos\theta_{N-1}^1\dots\sin\theta_1,\dots,\sin\theta_{N-1}^1\right),$$
$$(v_1)_{x_ix_i}=w^{a_1}_{\rho\rho}(\rho_1)\left(\frac{x_i-O^1_{x_i}}{\rho_1}\right)^2
+\frac{w^{a_1}_{\rho}(\rho_1)}{\rho_1}\left(1-\left(\frac{x_i-O^1_{x_i}}{\rho_1}\right)^2
-\frac{\partial O^1_{x_i}}{\partial x_i}\right)=$$
$$w^{a_1}_{\rho\rho}(\rho_1)\cos^2\theta_{N-1}\dots\sin^2\theta_{i-1}
+\frac{w^{a_1}_{\rho}(\rho_1)}{\rho_1}\left(1-\cos^2\theta_{N-1}\dots\sin^2\theta_{i-1}
-\frac{\partial O^1_{x_i}}{\partial x_i}\right).$$
Then
$$F(v_1)=\sum_{i=1}^{N-1}(v_1)_{x_ix_i}+(v_1)_{x_Nx_N}g((v_1)_{x_N}^2)=$$
$$w^{a_1}_{\rho\rho}(\rho_1)\left(\cos^2\theta_{N-1}^1+\sin^2\theta_{N-1}^1g\left((w^{a_1}_{\rho}(\rho_1))^2\sin^2\theta_{N-1}^1\right)\right)+$$
$$\frac{w^{a_1}_{\rho}(\rho_1)}{\rho_1}\left(N-2+\sin^2\theta_{N-1}^1
+\cos^2\theta_{N-1}^1g\left((w^{a_1}_{\rho}(\rho_1))^2\sin^2\theta_{N-1}^1\right)\right)-$$
$$\frac{w^{a_1}_{\rho}(\rho_1)}{\rho_1}\sum_{i=1}^{N-1}\frac{\partial O^1_{x_i}}{\partial x_i}\geq0$$
since $w^{a_1}(\rho_1)$ verifies equation (\ref{eqa}).

The sets $\mathcal{D}_0$ and $\mathcal{D}_1$ intersect on $\sin\theta_2(x_1,\dots,x_N)=\sin\theta_2^1(x_1,\dots,x_N)=\sqrt{a_1}$. For a point $(x_1,\dots,x_N)$ in this intersection we have
$$\rho_1(x_1,\dots,x_N)=\rho(x_1,\dots,x_N)-(R(1)-R(a_1)).$$
Hence, on $\mathcal{D}_0\cap\mathcal{D}_1$
$$R(a_1)\leq\rho_1(x_1,\dots,x_N)\leq R(a_1)+r$$
if and only if
$$R(1)\leq\rho(x_1,\dots,x_N)\leq R(1)+r,$$
and the functions $v_0$ and $v_1$ coincide:
$$v_1(x_1,\dots,x_N)=w^{a_1}(\rho_1(x_1,\dots,x_N))=w(\rho_1(x_1,\dots,x_N)+R(1)-R(a_1))=$$
$$w(\rho(x_1,\dots,x_N))=v_0(x_1,\dots,x_N).$$
The formula
$$\bar{v}(x_1,\dots,x_N)=\left\{\begin{array}{ll}
        v_0(x_1,\dots,x_N) & \mbox{ on }\mathcal{D}_0\\
        v_1(x_1,\dots,x_N) & \mbox{ on }\mathcal{D}_1\\
        \end{array}
\right.$$
defines a function $\bar{v}$ in $C^0(int(\mathcal{D}_0\cup\mathcal{D}_1))$. We claim that it is also in $C^1(int(\mathcal{D}_0\cup\mathcal{D}_1))$.

In fact, we have
$$\nabla v_0(x_1,\dots,x_N)=\frac{w_{\rho}(\rho(x_1,\dots,x_N))}{\rho(x_1,\dots,x_N)}(x_1,\dots,x_N),$$
$$\nabla v_1(x_1,\dots,x_N)=
\frac{w^{a_1}_{\rho}(\rho_1(x_1,\dots,x_N))}{\rho_1(x_1,\dots,x_N)}\left(x_1-O^1_{x_1},\dots,x_N-O^1_{x_N}\right)$$
and, on $\mathcal{D}_0\cap\mathcal{D}_1$,
$$\frac{1}{\rho_1(x_1,\dots,x_N)}\left(x_1-O^1_{x_1},\dots,x_N-O^1_{x_N}\right)=$$
$$\frac{1}{\rho(x_1,\dots,x_N)}(x_1,\dots,x_N).$$
On $int(\mathcal{D}_0)$ and $int(\mathcal{D}_1)$, the function $\bar{v}$ is of class $C^2$ and satisfies, pointwise, the inequality $F(\bar{v})\geq0$.
We claim that $\bar{v}$ is also in $W^{1,2}(int(\mathcal{D}_0\cup\mathcal{D}_1))$ and that it is a weak solution to $F(\bar{v})\geq0$ on $int(\mathcal{D}_0\cup\mathcal{D}_1)$.
In fact, for every $\eta\in C^\infty_0(int(\mathcal{D}_0\cup\mathcal{D}_1))$, applying the divergence theorem separately to $int(\mathcal{D}_0)$ and to $int(\mathcal{D}_1)$, we obtain
$$\int_{int(\mathcal{D}_0\cup\mathcal{D}_1)}
{\left[\mbox{div}\nabla_{\nabla v} L(\nabla\bar{v}(x))\eta(x)+\langle\nabla L(\nabla\bar{v}(x)),\nabla\eta(x)\rangle\right]dx}=$$
$$\int_{int(\mathcal{D}_0)\cup int(\mathcal{D}_1)}
{\left[\mbox{div}\nabla_{\nabla v} L(\nabla\bar{v}(x))\eta(x)+\langle\nabla L(\nabla\bar{v}(x)),\nabla\eta(x)\rangle\right]dx}=$$
$$\int_{\partial(int(\mathcal{D}_0))}{\eta(x)\langle\nabla L(\nabla\bar{v}(x)),\mbox{\bf n}(x)\rangle dl}+\int_{\partial(int(\mathcal{D}_1))}{\eta(x)\langle\nabla L(\nabla\bar{v}(x)),\mbox{\bf{n}}(x)\rangle dl}=$$
$$\int_{\mathcal{D}_0\cap\{\sin\theta_{N-1}=\sqrt{a_1}\}}{\eta(x)\langle\nabla L(\nabla\bar{v}(x)),\mbox{\bf n}(x)\rangle dl}+$$
$$\int_{\mathcal{D}_1\cap\{\sin\theta^1_{N-1}=\sqrt{a_1}\}}{\eta(x)\langle\nabla L(\nabla\bar{v}(x)),\mbox{\bf n}(x)\rangle dl}.$$
The last term equals zero, since $\bar{v}\in C^1(int(\mathcal{D}_0\cup\mathcal{D}_1))$. Hence, when if $\eta\geq0$, we have that
$$\int_{int(\mathcal{D}_0\cup\mathcal{D}_1)}{\langle\nabla L(\nabla\bar{v}(x)),\nabla\eta(x)\rangle dx}\leq0,$$
as we wanted to show.

Assuming defined $O^{n-2}(\theta_1,\dots,\theta_{N-2})$ and a function $v\in C^1(int(\mathcal{D}_0\cup\dots\cup\mathcal{D}_{n-2}))$, consider
$$O^{n-1}(\theta_1,\dots,\theta_{N-2})=(O^{n-1}_{x_1},\dots,O^{n-1}_{x_N})=$$
$$O^{n-2}(\theta_1,\dots,\theta_{N-2})+(R(a_{n-2})-R(a_{n-1}))$$
$$\left(\sqrt{1-a_{n-1}}\cos\theta_{N-2}\dots\cos\theta_1,\sqrt{1-a_{n-1}}\cos\theta_{N-2}\dots\sin\theta_1,\dots,\sqrt{a_{n-1}}\right).$$
Set
$$\rho_{n-1}(x_1,\dots,x_N)=\sqrt{(x_1-O^{n-1}_{x_1})^2+\dots+(x_N-O^{n-1}_{x_N})^2},$$
$$\sin\theta_{N-1}^{n-1}(x_1,\dots,x_N)=\frac{x_N-O^{n-1}_{x_N}}{\rho_{n-1}(x_1,\dots,x_N)}$$
$$\mathcal{D}_{n-1}=\left\{(x_1,\dots,x_N):R(a_n)<\rho_{n-1}(x_1,\dots,x_N)<R(a_{n-1})+r,\right.$$
$$\left.0=\sqrt{a_n}\leq\sin\theta_{N-1}^{n-1}(x_1,\dots,x_N)\leq\sqrt{a_{n-1}}\right\}$$
and define on $\mathcal{D}_{n-1}$ the function
$$v_{n-1}(x_1,\dots,x_N)=w^{a_{n-1}}(\rho_{n-1}(x_1,\dots,x_N)).$$
Set $\mathcal{D}=int(\mathcal{D}_0\cup\dots\cup\mathcal{D}_{n-1})$, the same considerations as before imply that the function
$$\bar{v}(x_1,\dots,x_N)=\left\{\begin{array}{ll}
        v_0(x_1,\dots,x_N) & \mbox{ on }\mathcal{D}_0\\
        \dots & \dots \\
        v_{n-1}(x_1,\dots,x_N) & \mbox{ on }\mathcal{D}_{n-1}\\
        \end{array}
\right.$$
is such that $\bar{v}\in W^{1,2}(\mathcal{D})\cap C^1(\mathcal{D})\cap C(\overline{\mathcal{D}})$ and it is a weak solution to $F(\bar{v})\geq0$ on $\mathcal{D}$. This completes the construction of $\bar{v}$ as a weak solution to $F(\bar{v})\geq0$ on $\mathcal{D}_0\cup\dots\cup\mathcal{D}_{N-1}$.

Set $O^*=(O^*_{x_1},\dots,O^*_{x_N})=(0,\dots,0,R(1)-l_N)$. We have that
$$\mathcal{D}_0\cup\dots\cup\mathcal{D}_{n-1}\subset\left\{(x_1,\dots,x_N):0\leq x_i\leq l+r
\mbox{ for }i=1,\dots,N-1,\right.$$
$$\left.O^*_{x_N}\leq x_N\leq O^*_{x_N}+l_N+r)\right\}.$$

\begin{figure}
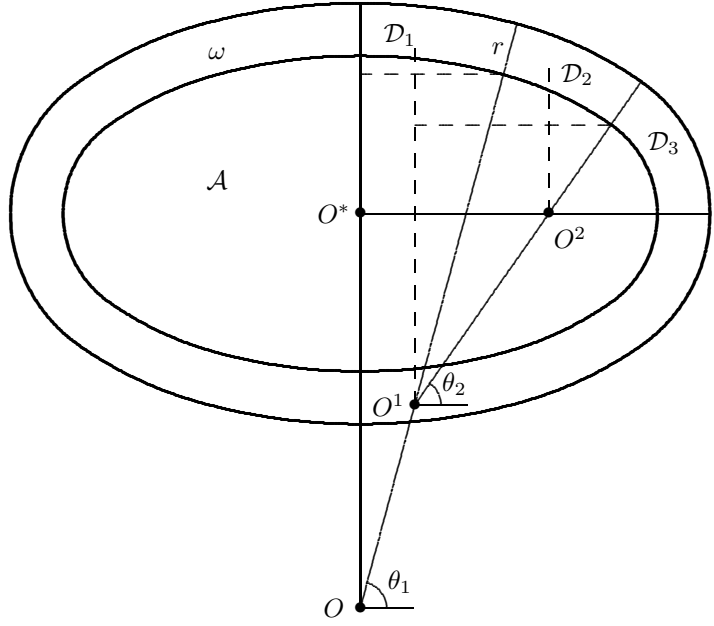


\begin{center}
\hspace{0cm}
\beginpicture

\setcoordinatesystem units < .35cm, .35cm>
\setplotarea x from -15 to  15, y from -16 to  10
\put {$\bullet$} at  0  0
\putrule from   0  8 to  0 -15
\plot   0 -15     5.952838  7.216293 /
\put {$\bullet$} at  0 -15
\plot   2.070553 -7.272593     10.6742  5.014687 /
\put {$\bullet$} at  2.070553 -7.272593
\putrule from   0  0 to  13.28468  0
\put {$\bullet$} at  7.162878  0
\setdashes
\plot   5.4352  5.284443     0  5.284443   /
\putrule from   2.070553 -7.272593 to  2.070553  6.727407
\putrule from   7.162878  0 to  7.162878  6
\plot   9.527046  3.376383     2.070553  3.376383   /
\setsolid
\putrule from   0 -15 to  2 -15
\putrule from   2.070553 -7.272593 to  4.070553 -7.272593
\circulararc -74.99997 degrees from  .258819 -14.03407 center at  0 -15
\circulararc -54.99998 degrees from  2.644129 -6.453442 center at  2.070553 -7.272593
\setplotsymbol({\large .})
\circulararc -30 degrees from -5.4352  5.284443 center at  0 -15
\circulararc -30 degrees from -5.952838  7.216293 center at  0 -15
\circulararc  30 degrees from -5.435199 -5.284443 center at  0  15
\circulararc  30 degrees from -5.952837 -7.216293 center at  0  15
\circulararc -20 degrees from  5.4352  5.284442 center at  2.070553 -7.272593
\circulararc -20 degrees from  5.952839  7.216294 center at  2.070553 -7.272593
\circulararc  20.00001 degrees from -5.4352  5.284442 center at -2.070553 -7.272593
\circulararc  20.00001 degrees from -5.952838  7.216294 center at -2.070553 -7.272593
\circulararc -20 degrees from -5.435199 -5.284443 center at -2.070553  7.272593
\circulararc -20 degrees from -5.952837 -7.216294 center at -2.070553  7.272593
\circulararc  20 degrees from  5.4352 -5.284442 center at  2.070553  7.272593
\circulararc  19.99999 degrees from  5.952838 -7.216294 center at  2.070553  7.272593
\circulararc -110 degrees from  9.527046  3.376383 center at  7.162878  0
\circulararc -110 degrees from  10.6742  5.014687 center at  7.162878  0
\circulararc -110 degrees from -9.527045 -3.376384 center at -7.162878  0
\circulararc -110 degrees from -10.6742 -5.014688 center at -7.162878  0
\setplotsymbol({\tiny .})
\put {$O^*$} at -1  0
\put {$O$} at -1 -15
\put {$O^1$} at  1.070553 -7.272593
\put {$O^2$} at  7.962878 -.9
\put {$\theta_1$} at  1.5 -14
\put {$\theta_2$} at  3.570553 -6.472593
\put {$r$} at  5.2352  6.284443
\put {$\omega$} at -5.4352  6.084442
\put {$\mathcal A$} at -5.4352  1.284442
\put {$\mathcal{D}_1$} at  1.5  6.8
\put {$\mathcal{D}_2$} at  8.2352  5.284443
\put {$\mathcal{D}_3$} at  11.52705  2.676383

\endpicture
\caption{The sets $\omega$ and $\mathcal A$ in the case $N=2$ and $n=3$.}
\label{fA}
\end{center}

\end{figure}

Define the full domain $\omega$ and the solution by symmetry with respect to the point $O^*$. Figure \ref{fA} shows this construction in dimension $N=2$ and for $n-1=2$. Hence the solution will be in $W^{1,2}(\omega)\cap C^1(\omega)\cap C(\overline{\omega})$ and a weak solution of $F(v)\geq0$ on $\omega$.

7) The previous construction yields a region $\mathcal{A}$ centered in $O^*$, a corresponding region $\omega$ and a function $v$ that solves (\ref{eqnN}). The change of coordinates $\hat x_1=x_1,\dots,\hat x_{N-1}=x_{N-1}$, $\hat x_N=x_N-O^*_N$, centers $\mathcal{A}$ at the origin and proves the theorem.\\

\endproof

In order to prove Theorem \ref{thsmp}, we need this further lemma.

\begin{lemma}\label{dist}

Consider the sets $\mathcal A$ and $\mathcal{R}(O^*,l,l_N)$, where $\mathcal A$, $O^*$, $l$, $l_N$ have been defined in Theorem \ref{thsmp}. Then, for every $p\in\partial\mathcal R$,
$$d(p,\overline{\mathcal A})<l_N.$$

\end{lemma}

\proof

\begin{figure}
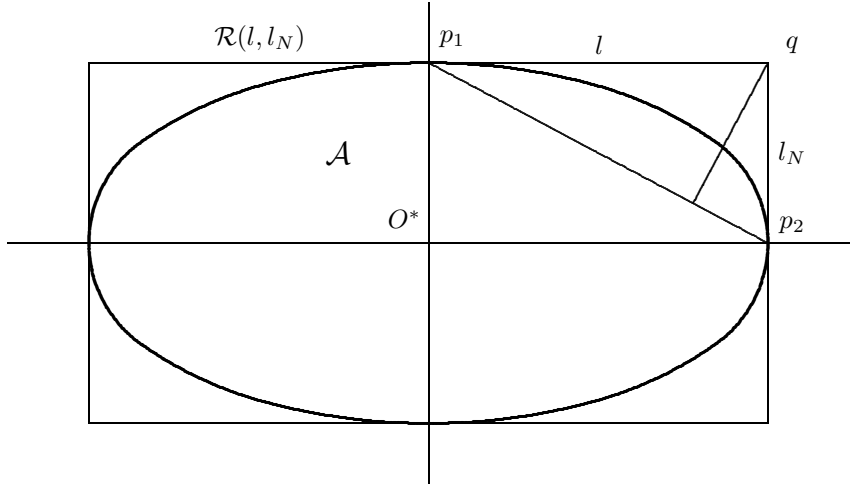

\begin{center}
\beginpicture

\setcoordinatesystem units < .4cm, .4cm>
\setplotarea x from -16 to  15, y from -8 to  8
\putrule from   14  0 to -14  0
\putrule from   0  8 to  0 -8
\setplotsymbol({\large .})
\circulararc -30 degrees from -5.4352  5.284443 center at  0 -15
\circulararc  30 degrees from -5.435199 -5.284443 center at  0  15
\circulararc -20 degrees from  5.4352  5.284442 center at  2.070553 -7.272593
\circulararc  20.00001 degrees from -5.4352  5.284442 center at -2.070553 -7.272593
\circulararc -20 degrees from -5.435199 -5.284443 center at -2.070553  7.272593
\circulararc  20 degrees from  5.4352 -5.284442 center at  2.070553  7.272593
\circulararc -110 degrees from  9.527046  3.376383 center at  7.162878  0
\circulararc -110 degrees from -9.527045 -3.376384 center at -7.162878  0
\setplotsymbol({\tiny .})
\putrule from   11.28468  6 to -11.28468  6
\putrule from  -11.28468  6 to -11.28468 -6
\putrule from  -11.28468 -6 to  11.28468 -6
\putrule from   11.28468 -6 to  11.28468  6
\plot   11.28468  0     0  6 /
\plot   11.28468  6     8.797606  1.322362   /
\put {\Large{$\mathcal A$}} at -3  3
\put {$\mathcal R(l,l_N)$} at -5.6  6.8
\put {$O^*$} at -.8  .8
\put {$q$} at  12.08468  6.6
\put {$p_1$} at  .8  6.8
\put {$p_2$} at  12.08468  .6
\put {$l$} at  5.68468  6.6
\put {$l_N$} at  12.08468  3

\endpicture
\caption{$\mathcal{A}$ and $\mathcal{R}(O^*,l,l_N)$ in the case $N=2$.}
\label{fdist}
\end{center}
\end{figure}

Set $q=O^*+(l,\dots,l,l_N)$. We prove that
$$d(q,\overline{\mathcal A})<l_N.$$
Set $p_i=O^*+(0,\dots,0,l,0,\dots,0)$ and let $\Pi^{N-1}$ the hyperplane passing through $p_1,\dots,p_N$. Since $\overline{\mathcal A}$ is convex and $p_i\in\overline{\mathcal A}$, we obtain that
$$d(q,\overline{\mathcal A})<d(q,\Pi^{N-1})<l_N.$$
See Figure \ref{fdist}.

\endproof

\proof[Proof of Theorem \ref{thsmp}]

a) Suppose that $u$ attains its minimum in $\Omega$, and assume $\min_\Omega u=0$ and set $\mathcal{C}=\{x\in\Omega:u(x)=0\}$. By contradiction, suppose that the open set $\Omega\setminus\mathcal{C}\neq\emptyset$.\\

b) Since $\Omega$ is a connected set, there exist $s\in\mathcal{C}$ and $R>0$ such that $B(s,R)\subset\Omega$ and $B(s,R)\cap(\Omega\setminus\mathcal{C})\neq\emptyset$.
Let $p\in B(s,R)\cap(\Omega\setminus\mathcal{C})$.
Consider the line $\overline{ps}$. Moving $p$ along this line, we can assume that $B(p,d(p,\mathcal{C}))\subset(\Omega\setminus\mathcal{C})$, and that there exists one point $z\in\mathcal{C}$ such that $d(p,\mathcal{C})=d(p,z)$.\\

c) Fix $r$:
$$0<r<\frac{d(p,\mathcal{C})}{32(N-1)K^2+\frac{7}{8}}.$$
Set
$$\epsilon(r)=\min\left\lbrace u(z):z\in\overline{B\left(p,d(p,\mathcal{C})-\dfrac{r}{4}\right)}\right\rbrace,$$
we have that $\epsilon(r)>0$, and we set $\epsilon=\min\{\epsilon(r),r\bar{\xi}\}$.\\

d) For $r$ and $\epsilon$ as defined in c), consider: $l$, $l_N$, $\mathcal{A}$ and $v$ as defined in Theorem \ref{subsol}.
Without loss of generality, since the set $\mathcal{A}$ is symmetric with respect to both coordinate axis, we can suppose that $\overline{pz}$ belongs to the first quadrant, i.e. that, for every $i=1,\dots,N$, $z_i\geq p_i$, where $z=(z_1,\dots,z_N)$ and $p=(p_1,\dots,p_N)$.

\begin{center}
\begin{figure}
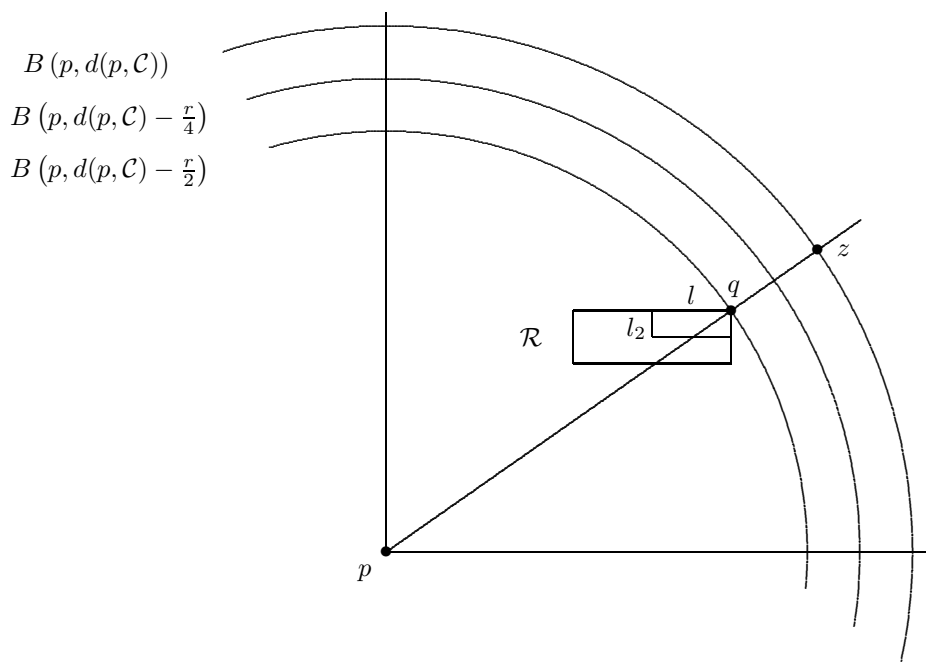

\beginpicture

\setcoordinatesystem units < .35cm, .35cm>
\setplotarea x from -12 to  21, y from -4 to  21
\putrule from   0  20.5 to  0  0
\putrule from   0  0 to  20.5  0
\circulararc -120 degrees from -6.180338  19.02113 center at  0  0
\circulararc -116 degrees from -5.262691  17.21349 center at  0  0
\circulararc -111 degrees from -4.410198  15.38019 center at  0  0
\setplotsymbol({\tiny .})
\setsolid
\putrule from   13.10643  9.177223 to  7.106433  9.177223
\putrule from   7.106433  9.177223 to  7.106433  7.177223
\putrule from   7.106433  7.177223 to  13.10643  7.177223
\putrule from   13.10643  7.177223 to  13.10643  9.177223
\plot   0  0     18.02135  12.61868 /
\setsolid
\putrule from   10.10643  9.177223 to  10.10643  8.177223
\putrule from   10.10643  8.177223 to  13.10643  8.177223
\setsolid
\put {$\bullet$} at  0  0
\put {$\bullet$} at  16.38304  11.47153
\put {$\bullet$} at  13.10643  9.177223
\put {$p$} at -.8 -.8
\put {$z$} at  17.38304  11.47153
\put {$B\left(p,d(p,\mathcal{C})-\frac{r}{2}\right)$} at  -10.5 14.5
\put {$B\left(p,d(p,\mathcal{C})-\frac{r}{4}\right)$} at  -10.5 16.5
\put {$B\left(p,d(p,\mathcal{C})\right)$} at  -11 18.5
\put {$q$} at  13.2 10
\put {$l$} at  11.6  9.8
\put {$l_2$} at  9.5  8.5
\put {$\mathcal R$} at  5.5  8.2

\endpicture
\caption{The set $\mathcal R=\mathcal{R}(q^*,l,l_2)$ in the case $N=2$.}
\label{fmax}

\end{figure}
\end{center}

Define the point $q$ on the segment $\overline{pz}$ such that $d(q,p)=d(p,z)-\dfrac{r}{2}$.
Set $q^*=q-(l,\dots,l,l_N)$, $\mathcal{R}(q^*,l,l_N)=q^*+\mathcal{R}(l,l_N)$, $\mathcal{A}^*=q^*+\mathcal{A}$ and $v^*(x+q^*)=v(x)$.
We first claim that
$$\mathcal{R}(q^*,l,l_N)\subset B\left(p,d(p,\mathcal{C})-\dfrac{r}{4}\right).$$
Let $t\in\mathcal{R}(q^*,l,l_N)$, then $t$ can be written as $(q_1-2\alpha_1l,\dots,q_{N-1}-2\alpha_{N-1}l,q_N-2\alpha_Nl_N)$, with $0\leq\alpha_i\leq1$, for $i=1,\dots,N$. Since $r<\frac{d(p,\mathcal{C})}{32(N-1)K^2+\frac{7}{8}}$, we have that
$$d(t,p)^2=\sum_{i=1}^N(q_i-2\alpha_il_i-p_i)^2=d(q,p)^2+\sum_{i=1}^N4\alpha_i^2l_i^2-\sum_{i=1}^N4\alpha_il_i(q_i-p_i)\leq$$
$$\left(d(p,\mathcal{C})-\dfrac{r}{2}\right)^2+\sum_{i=1}^N4l_i^2\leq\left(d(p,\mathcal{C})-\dfrac{r}{2}\right)^2+16(N-1)K^2r^2+\frac{r^2}{4}<\left(d(p,\mathcal{C})-\dfrac{r}{4}\right)^2.$$
See Figure \ref{fmax}.

Since $\mathcal{A}^*\subset\mathcal{R}(q^*,l,l_N)$, we have obtained that
$$\mathcal{A}^*\subset B\left(p,d(p,\mathcal{C})-\dfrac{r}{4}\right),$$
so that $u\geq\epsilon$ in $\partial\mathcal A^*$.

By Lemma \ref{dist},
$$d(q,\overline{\mathcal A^*})<l_N\leq\frac{r}{4},$$
we have that
$$d(z,\overline{\mathcal A^*})\leq d(z,q)+d(q,\overline{\mathcal A^*})<\frac{3}{4}r,$$
so that
$$z\in\omega^*=B(\mathcal A^*,r)\setminus\overline{\mathcal A^*}.$$

e) The function $v^*$ satisfies
$$\left\{\begin{array}{ll}
        v^*\mbox{ is a weak solution to }F(v)\geq0 & \mbox{ in }\omega^*\\
        v^*>0 & \mbox{ in }\omega^*\\
        v^*=0 & \mbox{ in }\partial B(\mathcal{A}^*,r)\\
        v^*\leq\epsilon & \mbox{ in }\partial\mathcal{A}^*.\\
        \end{array}
\right.
$$
Since $u,v^*\in W^{1,2}(\omega^*)\cap C(\overline{\omega^*})$, $v^*$ is a weak subsolution and $u$ is a weak solution to $F(u)=0$, and $v^*_{|\partial \omega^*}\leq u_{|\partial\omega^*}$, applying Lemma \ref{comp}, we obtain that $u\geq v^*$ in $\omega^*$. But $u(z)=0<v^*(z)$, a contradiction.\\

\endproof

\end{document}